\documentclass[11pt]{article}
\usepackage{latexsym,amssymb,amsmath}
\usepackage{amssymb,amsthm,upref,amscd}
\usepackage{color}
\usepackage[T1]{fontenc}
\begin{document}
\baselineskip=15pt

\numberwithin{equation}{section}

\pretolerance1000
\newtheorem{theorem}{Theorem}[section]
\newtheorem{lemma}[theorem]{Lemma}
\newtheorem{proposition}[theorem]{Proposition}
\newtheorem{corollary}[theorem]{Corollary}
\newtheorem{remark}[theorem]{Remark}
\newtheorem{definition}[theorem]{Definition}
\renewcommand{\theequation}{\thesection.\arabic{equation}}

\newcommand{\A}{\mathbb{A}}
\newcommand{\B}{\mathbb{B}}
\newcommand{\C}{\mathbb{C}}
\newcommand{\D}{\mathbb{D}}
\newcommand{\E}{\mathbb{E}}
\newcommand{\F}{\mathbb{F}}
\newcommand{\G}{\mathbb{G}}
\newcommand{\I}{\mathbb{I}}
\newcommand{\J}{\mathbb{J}}
\newcommand{\K}{\mathbb{K}}
\newcommand{\M}{\mathbb{M}}
\newcommand{\N}{\mathbb{N}}
\newcommand{\Q}{\mathbb{Q}}
\newcommand{\R}{\mathbb{R}}
\newcommand{\T}{\mathbb{T}}
\newcommand{\U}{\mathbb{U}}
\newcommand{\V}{\mathbb{V}}
\newcommand{\W}{\mathbb{W}}
\newcommand{\X}{\mathbb{X}}
\newcommand{\Y}{\mathbb{Y}}
\newcommand{\Z}{\mathbb{Z}}
\newcommand\ca{\mathcal{A}}
\newcommand\cb{\mathcal{B}}
\newcommand\cc{\mathcal{C}}
\newcommand\cd{\mathcal{D}}
\newcommand\ce{\mathcal{E}}
\newcommand\cf{\mathcal{F}}
\newcommand\cg{\mathcal{G}}
\newcommand\ch{\mathcal{H}}
\newcommand\ci{\mathcal{I}}
\newcommand\cj{\mathcal{J}}
\newcommand\ck{\mathcal{K}}
\newcommand\cl{\mathcal{L}}
\newcommand\cm{\mathcal{M}}
\newcommand\cn{\mathcal{N}}
\newcommand\co{\mathcal{O}}
\newcommand\cp{\mathcal{P}}
\newcommand\cq{\mathcal{Q}}
\newcommand\rr{\mathcal{R}}
\newcommand\cs{\mathcal{S}}
\newcommand\ct{\mathcal{T}}
\newcommand\cu{\mathcal{U}}
\newcommand\cv{\mathcal{V}}
\newcommand\cw{\mathcal{W}}
\newcommand\cx{\mathcal{X}}
\newcommand\ocd{\overline{\cd}}

\def\c{\centerline}
\def\ov{\overline}
\def\emp {\emptyset}
\def\pa {\partial}
\def\bl{\setminus}
\def\op{\oplus}
\def\sbt{\subset}
\def\un{\underline}
\def\al {\alpha}
\def\bt {\beta}
\def\de {\delta}
\def\Ga {\Gamma}
\def\ga {\gamma}
\def\lm {\lambda}
\def\Lam {\Lambda}
\def\om {\omega}
\def\Om {\Omega}
\def\sa {\sigma}
\def\vr {\varepsilon}
\def\va {\varphi}

\title{\bf Multi-bump solutions for a  Kirchhoff problem type }

\author{Claudianor O. Alves$^a$\thanks{C. O. Alves was partially supported by CNPq/Brazil
 301807/2013-2 and INCT-MAT, coalves@dme.ufcg.edu.br}\, , \, \   Giovany M. Figueiredo$^b$\thanks{Giovany M. Figueiredo was supported by 
CNPQ 302933/2014-0}\vspace{2mm}
\and {\small $a.$ Universidade Federal de Campina Grande} \\ {\small Unidade Acad\^emica de Matem\'{a}tica} \\ {\small CEP: 58429-900, Campina Grande - Pb, Brazil}\\
{\small $b.$ Universidade Federal do Par\'a} \\ { Faculdade de Matem\'atica,} \\ {\small CEP 66075-110, Bel\'em -Pa - Brazil.}}

\date{}
\maketitle

\begin{abstract}
In this paper, we are going to study the existence of solution for the following Kirchhoff problem
$$
\left\{ \begin{array}{l}
 M\biggl(\displaystyle\int_{\mathbb{R}^{3}}|\nabla u|^{2} dx +\displaystyle\int_{\mathbb{R}^{3}} \lambda a(x)+1)u^{2} dx\biggl) 
\biggl(- \Delta u + (\lambda a(x)+1)u\biggl)  = f(u)  \mbox{ in } \,\,\, \mathbb{R}^{3}, \\
\mbox{}\\
u \in H^{1}(\mathbb{R}^{3}).
\end{array}
\right.
$$
Assuming that the nonnegative function $a(x)$ has a potential well with  $int (a^{-1}(\{0\}))$  consisting of 
$k$ disjoint components $\Omega_1, \Omega_2, ....., \Omega_k$ and the nonlinearity $f(t)$ has a subcritical growth, 
we are able to establish the existence of positive multi-bump solutions by variational methods.

 \vspace{0.3cm}


\vspace{0.3cm}

\end{abstract}

\section{Introduction}

In the present paper, we are interested in showing the existence of multi-bump solution for the following class of Kirchhoff problem
$$
\left\{ \begin{array}{l}
 M\biggl(\displaystyle\int_{\mathbb{R}^{3}}|\nabla u|^{2} dx +\displaystyle\int_{\mathbb{R}^{3}} (\lambda a(x)+1)u^{2} dx\biggl) 
\biggl(- \Delta u + (\lambda a(x)+1)u\biggl)  = f(u)  \mbox{ in } \,\,\, \mathbb{R}^{3}, \\
\mbox{\hspace{13,7 cm}}{(P)_\lambda}\\
u \in H^{1}(\mathbb{R}^{3}),
\end{array}
\right.
$$
where $\lambda >0 $ is a positive parameter and $M, a $ and $f$ are functions verifying some conditions, which will be fixed below. 

The  function  $M: \mathbb{R}_{+} \to \mathbb{R}_{+}$ belongs to $C^{1}(\mathbb{R}, \mathbb{R})$ and satisfies the following conditions:

\begin{description}

\item[($M_{1}$)] The function $M$ is increasing and $0 < M(0)=: m_{0}$.

\item[($M_{2}$)] The function $t\mapsto\displaystyle\frac{M(t)}{t}$ is decreasing.

\end{description}

A typical example of function verifying the assumptions
$(M_{1})-(M_{2})$ is given by
$$
\displaystyle  M(t)=m_{0}+bt, \ \ \mbox{where} \ \ m_{0}>0  \ \ \mbox{and} \ \
b>0.
$$
This is the example that was considered in \cite{kirchhoff}. More generally, each
function of the form
$$
\displaystyle  M(t)=m_{0}+bt+\displaystyle \sum_{i=1}^{k}b_{i}t^{\gamma_{i}}
$$
with $b_{i}\geq 0$ and $\gamma_{i}\in (0,1)$ for all $i\in
\{1,2,\ldots, k\}$ verifies the hypotheses $(M_{1})-(M_{2})$. An
another example is $M(t) = m_{0} + ln(1 + t)$.

Related to function $a(x)$, we assume the following conditions:

\vspace{0,3 cm}
\noindent $(a_0)$ \,\,\, $a(x) \geq 0, \quad \forall x \in \mathbb{R}^{N}$. \\
\noindent $(a_1)$ \,\,The set $int (a^{-1}(\{0\}))$ is nonempty and there are disjoint open components $\Omega_1, \Omega_2, ....., \Omega_k$ such that
\begin{equation} \label{a1}
int (a^{-1}(\{0\}))= \cup_{j=1}^{k} \Omega_j
\end{equation}
and
\begin{equation} \label{a2}
dist( \Omega_i, \Omega_j)>0 \,\,\, \mbox{for} \,\,\, i \not= j,\,\,\, i, j=1,2,\cdots, k .
\end{equation}

Finally, the function $f$ is a continuous function satisfying: 
\begin{enumerate}
\item[$(f_1)$\ ]
$\displaystyle\lim_{s\rightarrow 0 } \frac{f(s)}{s} = 0$,
\item[$(f_2)$\ ] $\displaystyle\lim_{|s|\rightarrow +\infty } \frac{f(s)}{s^5} = 0$,
\item[$(f_3)$\ ] There exists $\theta >4$ such that
\[
0 < \theta F(s) \leq sf(s)\quad \forall  s\in \mathbb{R}\setminus \{0\}.
\]
\item[$(f_4)$\ ] $\displaystyle \frac{f(s)}{s^3} $ is increasing in $s>0$ and decreasing in $s<0$.
\end{enumerate}

\vspace{0.5 cm}

Related to problem $(P)_\lambda$, we have the problem 
$$
\left \{ \begin{array}{l}
-M\left(\displaystyle\int_\Omega \mid\nabla u\mid^{2} dx\right)\Delta u =f(u) \ \mbox{in}\ \Omega,\\
u\in H^{1}_{0}(\Omega)                                    
\end{array}
\right.\leqno{(*)}
$$
where $\Omega \subset \mathbb{R}^{N}$ is a bounded domain. This type of problem is called {\it Kirchhoff problem}, because of the presence of the term
$M\left(\displaystyle\int_{\Omega}|\nabla u|^{2} dx \right)$. Indeed, this operator
appears in the Kirchhoff
equation \cite{kirchhoff}, which arises in nonlinear vibrations, namely
$$
\left\{
\begin{array}{l}
 u_{tt}-M\left(\displaystyle \int_{\Omega}|\nabla u|^{2} dx \right)\Delta u =
g(x,u) \ \mbox{in} \ \Omega \times (0,T)\\
u=0 \ \mbox{on} \ \partial\Omega \times (0,T)\\
u(x,0)=u_{0}(x)\ \ , \ \ u_{t}(x,0)=u_{1}(x).
\end{array}
\right.
$$

The reader may consult \cite{alvescorrea}, \cite{alvescorreama}, \cite{Cristian}, \cite{ma}
and the references therein, for more physical motivation on Kirchhoff problem. 

We would like point out that in the last years many authors have studied this type of problem in bounded or unbounded domains, see for example, \cite{AlvesFigueiredo}, \cite{Azzollini}, \cite{Cammaroto}, \cite{Cheng}, \cite{Chen}, \cite{SChen}, \cite{jmaa}, \cite{Cristian}, \cite{Giovany1}, \cite{Giovany2}, \cite{He}, \cite{Li}, \cite{Liao}, \cite{Liu}, \cite{Liu1}, \cite{Naimen}, \cite{Naimen1},  \cite{Xiang}, \cite{Wang1}, \cite{Wang2} and reference therein. For  solutions that change sign (nodal solution) we would like to cite \cite{Giovany3}, \cite{Mao}, \cite{Mao1}, \cite{Shuai} and \cite{Zhang}. 

The motivation to study the problem $(P)_\lambda$ comes from of a paper due to Ding and Tanaka \cite{DingTanaka}, which has studied  $(P)_\lambda$ assuming $M(t)=1$ and $f(t)=|t|^{q-1}t$. In that interesting paper, the authors considered the existence of positive multi-bump solution for the problem
\begin{equation} \label{LL2}
\left\{ \begin{array}{ll}
  - \Delta u + (\lambda a(x)+Z(x))u = u^{q} \,\,\,  \mbox{ in } \,\,\, \mathbb{R}^{N},\\
u \in H^{1}(\mathbb{R}^{N}), \\
\end{array}
\right.
\end{equation}
$ q \in (1, \frac{N+2}{N-2} ) $ if $ N \ge 3 $; $ q \in (1, \infty) $ if $ N = 1, 2 $. The authors showed that the above problem has at least $ 2^k-1 $ solutions $u_\lambda$ for large values of $ \lambda $. More precisely, for each non-empty subset $ \Upsilon $ of $ \{ 1,\ldots,k \} $, it was proved that, for any sequence $ \lambda_n \to \infty $ we can extract a subsequence $( \lambda_{n_i}) $ such that $( u_{ \lambda_{n_i} } )$ converges strongly in $ H^1 \big( \mathbb R^N \big) $ to a function $ u $, which satisfies $ u = 0 $ outside $ \Omega_\Upsilon = \bigcup_{ j \in \Upsilon } \Omega_j $ and $ u_{|_{\Omega_j}}, \, j \in \Upsilon $, is a least energy solution for
\begin{equation} \label{LL}
   \begin{cases}
		  - \Delta u + Z(x) u  = u^q, \text{ in } \Omega_j, \\
		  u \in H^1_0 \big( \Omega_j \big), \, u > 0, \text{ in } \Omega_j.
	 \end{cases}
\end{equation}

Involving the Kirchhoff problem with potential wells, there are not so many existing papers. As far as we know, the only paper that considered the existence of solutions for $(P)_{\lambda}$ is due to Liang and Shi \cite{SS}. Unfortunately, we believe that the    Section V  of the above paper has a mistake, which  commits the proof of the their main result, to be more precisely, we have observed that the numbers $c_j$ and $c_{\lambda,j}$ considered in that work are not a good choice for this class of problem, and also, the proof of Lemma 5.3 is not correct, because the authors have forgotten that the Kirchhoff problem has a nonlocal term involving the function $M$, which is very sensitive for some estimates. Motivated by \cite{DingTanaka} and \cite{SS}, we intend in the present paper to show how we can work with this nonlocal term to get  a positive multi-bump solution for $(P)_{\lambda}$.  Here, we will adapt an idea used  by  Alves and Yang \cite{AY} to show the  existence of multi-bump solution for the following Schr\"odinger-Poisson system
$$
\left\{ \begin{array}{ll}
  - \Delta u + (\lambda a(x)+1)u+ \phi u = f(u)  \mbox{ in } \,\,\, \mathbb{R}^{3},\\
-\Delta \phi=u^2  \mbox{ in } \,\,\, \mathbb{R}^{3}.\\
\end{array}
\right.
$$ 

\vspace{0.5 cm}

Our main result is the following

\begin{theorem} \label{main}
   Assume that $ (M_1), (M_2), (a_0),(a_1)$ and $ (f_1)-(f_4) $ hold. Then, there exist $ \lambda_0 > 0 $ 
	with the following property: for any non-empty subset $ \Upsilon $ of $ \{1, 2, . . . , k \} $ and 
	$ \lambda \ge \lambda_0 $,  problem $ \big( P_\lambda \big) $ has a positive solution $u_\lambda$. 
	Moreover, if we fix the subset $ \Upsilon $, then for any sequence $ \lambda_n \to \infty $ 
	we can extract a subsequence $ (\lambda_{n_i}) $ such that $ (u_{ \lambda_{n_i} }) $ 
	converges strongly in $ H^{1}(\mathbb R^3) $ to a function $ u $, 
	which satisfies $ u = 0 $ outside $ \Omega_\Upsilon = \cup_{ j \in \Upsilon } \Omega_j $, and $ u_{|_{ \Omega_\Upsilon}}$ 
	is a least energy solution for the nonlocal problem

$$
\left\{
\begin{array}{l}
M\biggl(\displaystyle\int_{\Omega_\Upsilon}|\nabla u|^{2} dx + \displaystyle\int_{\Omega_\Upsilon} u^{2} dx \biggl) 
\biggl(- \Delta u +u \biggl)= f(u)  \mbox{ in } \,\,\, \Omega_\Upsilon,\\
u(x)>0 \,\,\, \forall x \in \Omega_j \,\,\, \mbox{and} \,\,\, \forall j \in \Upsilon,\\
u \in H^{1}_{0}(\Omega_\Upsilon).
\end{array}
\right.
\eqno{(P)_{\infty, \Upsilon}}
$$
\end{theorem}

The paper is organized as follows. In the next section, we prove some technical lemmas
and the existence of least energy solution for $(P)_{\infty,\Upsilon}$. In Section 3, we study an auxiliary problem. A compactness result for the energy functional associated with the auxiliary problem is showed in Section 4. Some estimates involving the solutions of auxiliary problem are showed in Section 5, and in Section 6, we build a special minimax value for the functional energy associated to the  auxiliary problem. 

\section { The problem $(P)_{\infty,\Upsilon}$ }

In the sequel, let us denote by $\widehat{M}$ and ${F}$ the following functions
$$
\widehat{M}(t)=\displaystyle\int^{t}_{0}M(s) ds \quad \mbox{and} \quad  {F}(t)=\displaystyle\int^{t}_{0}f(s) ds.
$$ 

In the proof of Theorem \ref{main}, we need to study the existence of least energy solution for problem $(P)_{\infty,\Upsilon}$. The main idea is to prove that the energy functional $J$ associated with nonlocal problem $(P)_{\infty, \Upsilon}$ given by
$$
J(u)=\frac 12 \widehat{M}\biggl(\int_{\Omega_\Upsilon}(|\nabla u|^{2}+|u|^{2})dx\biggl) -\int_{\Omega_\Upsilon}F(u)dx,
$$
assumes a minimum value on the set
$$
\mathcal M_{\Upsilon}=\{u\in \mathcal N_{\Upsilon}: J'(u)u_j=0 \mbox { and }
u_{j}\neq 0 \,\,\, \forall j \in \Upsilon \}
$$
where $u_{j}=u_{|_{ \Omega_j}}$ and $\mathcal{N}_{\Upsilon}$ is the corresponding Nehari manifold defined by
$$
\mathcal{N}_{\Upsilon} = \{ u\in
H^1_0(\Omega_\Upsilon)\setminus\{0\}\, :\, J'(u)u=0\}.
$$
More precisely, we will prove that there is $w \in \mathcal{M}_\Upsilon$ such that
$$
J(w)=\inf_{u \in \mathcal{M}_\Upsilon}J(u).
$$
After, we use  the implicit function theorem to prove that $w$ is a
critical point of $J$, and so, $w$ is a least energy solution
for $(P)_{\infty,\Upsilon}$. The main feature of the least energy solution $w$ is that $w(x)>0 \,\,\, \forall x \in \Omega_j \,\,\, \mbox{and} \,\,\, \forall j \in \Upsilon$, which will be used to describe the existence of multi-bump solutions.

Since we intend to look for positive solutions, throughout this paper we assume that
$$
f(s)=0, \,\,\, s \leq 0.
$$

Moreover, notice that by $(f_{1})$ and $(f_{2})$,  given
$\epsilon>0$, there exists $C_{\epsilon}>0$ such that
\begin{eqnarray}\label{estimativaf1}
f(t)t \leq \epsilon|t|^{2} + C_{\epsilon}|t|^{6}.
\end{eqnarray}

In what follows, to show in details the idea of the existence of least energy solution for $(P)_{\infty,\Upsilon}$, we will consider $\Upsilon=\{1,2\}$. Moreover, we  will denote by $\Omega$, $\mathcal{N}$ and $ \mathcal M$ the sets $\Omega_\Upsilon$, $\mathcal{N}_\Upsilon$ and $ \mathcal{M}_\Upsilon $ respectively. Thereby,
$$
\Omega=\Omega_1 \cup \Omega_2,
$$
$$
\mathcal{N} = \{ u\in
H^1_0(\Omega)\setminus\{0\}\, :\, J'(u)u=0\}
$$
and
$$
\mathcal M=\{u\in \mathcal N: J'(u)u_1=J'(u)u_2=0 \mbox { and }
u_{1},u_{2}\neq 0 \},
$$
with $u_{j}=u_{|_{ \Omega_j}}$, $j=1,2.$

\subsection{Technical lemmas}

Hereafter, let us denote by $||\,\,\,||$, $||\,\,\,||_1$ and $||\,\,\,||_2$  the norms in $H^{1}_0(\Omega)$, $H^{1}_0(\Omega_1)$ and $H^{1}_0(\Omega_2)$ given by
$$
||u||=\left(\int_{\Omega}(|\nabla u|^{2}+|u|^{2})dx\right)^{\frac{1}{2}},
$$
$$
||u||_1=\left(\int_{\Omega_1}(|\nabla u|^{2}+|u|^{2})dx\right)^{\frac{1}{2}}
$$
and
$$
||u||_2=\left(\int_{\Omega_2}(|\nabla u|^{2}+|u|^{2})dx\right)^{\frac{1}{2}}
$$
respectively.

\vspace{0.5 cm}
In order to show that the set $\mathcal M$ is not empty, we need of the following Lemma.
\begin{lemma}\label{lema3}
Let  $v\in H_0^1(\Omega)$ with $v_{j}\neq 0$ for $j=1,2$. Then,
there are $t,s>0$ such that $J'(tv_1+sv_2)v_1=0$ and
$J'(tv_1+sv_2)v_2=0$.
\end{lemma}

\noindent {\bf{Proof.}} Let
$V:(0,+\infty)\times(0,+\infty)\rightarrow \mathbb{R}^{2}$ be a
continuous function given by
$$
V(t,s)=(J'(tv_{1}+sv_{2})(tv_{1}),J'(tv_{1}+sv_{2})(sv_{2})).
$$
Note that
\begin{eqnarray}\label{prausar1}
J'(tv_{1}+sv_{2})(tv_{1})=t^{2}M(t^{2}\|v_{1}\|^{2}_{1}+s^{2}\|v_{2}\|^{2}_{2})\|v_{1}\|^{2}_{1}-\displaystyle\int_{\Omega_{1}}f(tv_{1})tv_{1}
dx.
\end{eqnarray}

Using $(M_{1})$, (\ref{estimativaf1}) and Sobolev's embedding in
(\ref{prausar1}), we have
$$
J'(tv_{1}+sv_{2})(tv_{1})\geq (m_{0}-\epsilon
C)t^{2}\|v_{1}\|^{2}_{1}-t^{q}C_{\epsilon}C \|v_{1}\|^{6}_{1},
$$
for some $C>0$. Thus, there exists $r>0$ sufficiently small such
that
$$
J'(rv_{1}+sv_{2})(rv_{1})>0, \ \ \mbox{for all} \ \ s>0.
$$
The same idea yields 
$$
J'(tv_{1}+rv_{2})(rv_{2})>0, \ \ \mbox{for all} \ \ t>0.
$$
On the other hand, by $(M_{2})$, there exists $K_{1}>0$ such that
\begin{eqnarray}\label{estimativaporcima}
M(t)\leq M(1)t+K_{1}, \ \ \mbox{for all} \ \ t \geq 0
\end{eqnarray}
and by $(f_{3})$, there are $K_{2}, K_{3}>0$ such that
\begin{eqnarray}\label{estimativaf2}
F(t)\geq K_{2}t^{\theta}-K_{3}.
\end{eqnarray}
Using (\ref{estimativaporcima}), (\ref{estimativaf2}) and $(f_{3})$
in (\ref{prausar1}), we derive that
\begin{eqnarray*}
J'(tv_{1}+sv_{2})(tv_{1})&\leq & t^{4}M(1)\|v_{1}\|^{4}_{1}+
t^{2}s^{2}M(1)\|v_{1}\|^{2}_{1}\|v_{2}\|^{2}_{2}+K_{1}t^{2}\|v_{1}\|^{2}_{1}\\
&-&
\frac{t^{\theta}}{\theta}K_{2}\displaystyle\int_{\Omega}|v_{1}|^{\theta}_{1}
dx+K_{3}|\Omega_{1}|,
\end{eqnarray*}
where $|\Omega_{1}|$ denotes the Lebesgue measure of $\Omega_{1}$. Thus,
since $\theta>4$, for $R>0$ sufficiently large, we get
$$
J'(Rv_{1}+sv_{2})(Rv_{1})<0, \ \ \mbox{for all} \ \ s\leq R.
$$
Arguing of the same way, we also have 
$$
J'(tv_{1}+Rv_{2})(Rv_{2})<0, \ \ \mbox{for all} \ \ t\leq R.
$$
In particular,
$$
J'(rv_{1}+sv_{2})(rv_{1})>0 \ \ \mbox{and}  \ \
J'(tv_{1}+rv_{2})(rv_{2})>0, \ \ \mbox{for all} \ \ t,s \in [r,R]
$$
and
$$
J'(Rv_{1}+sv_{2})(Rv_{1})<0 \ \ \mbox{and}  \ \
J'(tv_{1}+Rv_{2})(Rv_{2})<0, \ \ \mbox{for all} \ \ t,s \in [r,R].
$$
Now the lemma follows applying Miranda's Theorem \cite{Miranda}.
\qed

\hspace{0.5cm}

As an immediate consequence of the last lemma, we have the following corollary

\begin{corollary} \label{CM}

The set $\mathcal M$ is not empty.

\end{corollary}

Next, we will show some technical lemmas.

\begin{lemma}\label{lema2}
There exists $\rho>0$ such that
\begin{enumerate}
\item [(i)] $J(u)\geq  \frac{(\theta -4)}{4\theta}m_{0} ||u||^2$ and $||u||\geq \rho, \forall u\in \mathcal N$;
\item [(ii)] $||w_j||_j\geq\rho, \,\, \forall w \in \mathcal M$ and $j=1,2$, where $w_j=w|_{\Omega_j}, j=1,2$.
\end{enumerate}
\end{lemma}

\noindent {\bf{Proof.}} From definition of $\widehat{M}$ and
$(M_{2})$, 
\begin{eqnarray}\label{estimativaM1}
\widehat{M}(t)\geq \frac{1}{2}M(t)t, \ \ \mbox{for all}\ \ t\geq 0.
\end{eqnarray}
Now, a simple computation together with (\ref{estimativaM1}) gives 
\begin{eqnarray}\label{estimativaM2}
\frac{1}{2}\widehat{M}(t)- \frac{1}{\theta}M(t)t \geq \frac{(\theta
-4)}{4\theta}m_{0}t, \ \ \mbox{for all}\ \ t\geq 0.
\end{eqnarray}
Thus, by $(f_{3})$ and (\ref{estimativaM2}), 
$$
J(u)= J(u)-\frac{1}{\theta}J'(u)u\geq \frac{(\theta
-4)}{4\theta}m_{0}\|u\|^{2}, \ \ \mbox{for all} \ \ u \in
\mathcal{N}. 
$$
Gathering definition of $\mathcal{N}$, $(M_{1})$,
(\ref{estimativaf1}) and Sobolev's embedding, it follows that
$$
0< \rho:=\biggl[\bigl(m_{0}-\frac{\epsilon
C_{2}}{C_{1}}\bigl)\frac{1}{C_{\epsilon}}\biggl]^{1/(q-2)}\leq
\|u\|,
$$
for all $u \in \mathcal{N}$ and for some $C_{1}, C_{1}>0$.

From $(M_{1})$,
$$
M(\|w_{j}\|^{2}_{j})\leq M(\|w\|^{2}), \quad \forall w \in \mathcal{M}.
$$
Thus, 
\begin{eqnarray}\label{beleza}
J'(w_{j})w_{j}\leq 0, \ \ \mbox{for all} \ \ w \in \mathcal{M},
\end{eqnarray}
implying that 
$$
0< \rho\leq \|w_{j}\|_{j}.
$$

 \qed

 \begin{lemma}\label{lema2x}
If $(w_n)$ is a bounded sequence in  $\mathcal M$ and $q\in
(2,6)$, we have
$$
\liminf_n \int_{\Omega_j}|w_{n,j}|^{p}dx>0 \,\,\, j=1,2.
$$
where $w_{n,j}=w_n|_{\Omega_j}$ for $j=1,2$. \end{lemma}

\noindent {\bf{Proof.}} 
Notice that by $(f_{1})$ and $(f_{2})$,  given
$\epsilon>0$, there exist $C>0$ and $q\in
(2,6)$ such that
\begin{eqnarray}\label{estimativaf101}
f(t)t \leq \epsilon|t|^{2} + C|t|^{q}+ \epsilon|t|^{6}.
\end{eqnarray}
Therefore,
$$
0<m_{0}\rho^{2}\leq M(\|w_{n,j}\|^{2}_{j})\|w_{n,j}\|^{2}_{j}\leq
\epsilon \displaystyle\int_{\Omega}|w_{n,j}|^{2}_{j} \ dx +
C \displaystyle\int_{\Omega}|w_{n,j}|^{q} \ dx + \epsilon \displaystyle\int_{\Omega}|w_{n,j}|^{6}_{j} \ dx.
$$
Since $(w_{n})$ is bounded, there is $\tilde{C}>0$ such that
$$
0<m_{0}\rho^{2}\leq \epsilon \tilde{C}+ C
\displaystyle\int_{\Omega}|w_{n,j}|^{q} \ dx.
$$
Now, the result follows choosing $\epsilon$ small enough. \qed

\subsection { Existence of least energy solution for $(P)_{\infty,\Upsilon}$  }

At this point, some useful remarks follow. First of all,  let us
observe that, from $(M_{2})$, 
\begin{eqnarray}\label{relacaoderivadafuncaoM}
M'(t)t\leq M(t), \ \ \mbox{for all} \ \ t\geq 0,
\end{eqnarray}
from where it follows that
\begin{eqnarray}\label{monotonicy1}
t\mapsto \frac{1}{2}\widehat{M}(t)- \frac{1}{4}M(t)t \ \ \mbox{is
increasing}.
\end{eqnarray}

Now, by $(f_{4})$, 
\begin{eqnarray}\label{relacaoderivadafuncaof}
f'(t)t\geq 3f(t),\ \ \mbox{for all} \ \ |t|\geq 0,
\end{eqnarray}
implying that
\begin{eqnarray}\label{monotonicy2}
t\mapsto \frac{1}{4}f(t)t- F(t) \ \ \mbox{is increasing, for all} \
\ |t|>0.
\end{eqnarray}

In this subsection, our main goal is to prove the following result:

\begin{theorem} \label{T2}
Assume that $(f_1)-(f_4)$ hold. Then problem $(P)_{\infty,\Upsilon}$ possesses a positive least energy solution on the set $\mathcal{M}$.
\end{theorem}

\noindent {\bf{Proof.}} We will prove the existence of $w \in \mathcal{M}$
in which the infimum of $J$ is attained on $ \mathcal{M}$. After, using the implicit function theorem,  
we show that $w$ is a critical point of $J$, from where it follows that $w$ is a least energy  solution of $(P)_{\infty,\Upsilon}$. 

\medskip

First of all, by Lemma \ref{lema2}, there
exists $c_{0} \in \mathbb{R}$ such that
$$
0<c_{0}=\displaystyle\inf_{v \in \mathcal{M}}J(v).
$$

\noindent Thus, by Corollary \ref{CM}, there exists a  minimizing sequence $(w_{n})$ in
$\mathcal{M}$, which  is bounded, by Lemma \ref{lema2}.
Hence, by Sobolev Imbedding Theorem, without loss of
generality, we can assume up to a subsequence that there exist $w 
 \in H^{1}_{0}(\Omega)$ such that
\begin{eqnarray*}
w_{n}\rightharpoonup w \ \  \mbox{in} \ \ H^{1}_{0}(\Omega), \ \  w_{n} \rightarrow w      \ \  \mbox{in} \quad L^q(\Omega) 
\ \   \mbox{with} \ \ q \in (1,6) \ \ \mbox{and} \ \  w_{n}(x) \rightarrow w(x) \ \ \mbox{a.e in} \ \ \Omega.
\end{eqnarray*}

Then, $(f_2)$ combined with the \emph{compactness lemma of Strauss}
\cite[Theorem A.I, p.338]{bl} gives
$$
\lim_n \int_{\Omega_j}|w_{n,j}|^{q}dx = \int_{\Omega_j}|w_j|^{q}dx,
$$
$$
\lim_n \int_{\Omega_j} w_{n,j}f(w_{n,j})dx=\int_{\Omega_j} w_j f(w_j)dx
$$
and
$$
\lim_n \int_{\Omega_j}  F(w_{n,j})dx=\int_{\Omega_j} F(w_j)dx,
$$
from where it follows together with Lemma \ref{lema2x} that $w_j\neq 0$ for $j=1,2$. Thereby, by Lemma \ref{lema3}, there are $t,s>0$ verifying
$$
J'(tw_1+sw_2)w_1=0 \,\,\, \mbox{and} \,\,\, J'(tw_1+sw_2)w_2=0.
$$

Now, let us prove that $t, s
\leq 1$. First of all, we observe that subcritical growth of $f$ loads to 
growth, we get
$$
\displaystyle\int_{\Omega}f(w_{n,j})w_{n,j} dx
\rightarrow \displaystyle\int_{\Omega}f(w_{j})w_{j} dx
$$
and
$$
\displaystyle\int_{\Omega}F(w_{n,j}) dx  \rightarrow
\displaystyle\int_{\Omega}F(w_{j}) dx.
$$

\medskip

\noindent Thus,  as $ J'(w_{n})w_{n,j} = 0$,  
$$
M(\|w_{n}\|^{2})\|w_{n,1}\|^{2} = \displaystyle\int_{\Omega_{1}}f(w_{n,1})w_{n,1},
$$
or equivalently,
\begin{eqnarray*}
\frac{M(\|w_{n}\|^{2})}{\|w_{n}\|^{2}}\|w_{n,1}\|^{2}\|w_{n}\|^{2} = \displaystyle\int_{\Omega_{1}}
\frac{f(w_{n,1})}{w_{n,1}^{3}}{w_{n,1}^{4}} dx.
\end{eqnarray*}
Taking the limit in the above equality, we find
\begin{eqnarray}\label{comparacao}
\frac{M(\|w\|^{2})}{\|w\|^{2}}\|w_{1}\|^{2}\|w\|^{2} \leq  \displaystyle\int_{\Omega_{1}}
\frac{f(w_{1})}{w_{1}^{3}}{w_{1}^{4}} dx.
\end{eqnarray}
On the other hand, as $ J'(tw_{1}+ sw_{2})tw_{1} = 0$, we must have
$$
M(\|tw_{1}+sw_{2}\|^{2})\|tw_{1}\|^{2} = \displaystyle\int_{\Omega_{1}}f(tw_{1})tw_{1} dx. 
$$
Without generality, we can suppose $s\leq t$. Hence,
\begin{eqnarray}\label{comparacao1}
\frac{M(t^{2}\|w\|^{2})}{t^{2}\|w\|^{2}}\|w_{1}\|^{2}\|w\|^{2} \geq  \displaystyle\int_{\Omega_{1}}
\frac{f(tw_{1})}{(tw_{1})^{3}}{w_{1}^{4}} dx.
\end{eqnarray}
Combining (\ref{comparacao}) with (\ref{comparacao1}), we derive

\begin{eqnarray*}
\biggl[\frac{M(t^{2}\|w\|^{2})}{t^{2}\|w\|^{2}}-\frac{M(\|w\|^{2})}{\|w\|^{2}}\biggl]\|w_{1}\|^{2}\|w\|^{2} \geq  
\displaystyle\int_{\Omega_{1}}
\biggl[\frac{f(tw_{1})}{(tw_{1})^{3}}-\frac{f(w_{1})}{(w_{1})^{3}}\biggl]{w_{1}^{4}} dx.
\end{eqnarray*}
Gathering $(M_{2})$ and $(f_{4})$, we ensure that $0<s\leq t\leq 1$.

In the next step, we will show that $J(tw_{1}+sw_{2})=c_{0}$. Since $tw_{1}+sw_{2} \in \mathcal{M}$ and $t, s \leq 1 $, from 
(\ref{monotonicy1}) and (\ref{monotonicy2}), 
\begin{eqnarray*}
c_{0}\leq J(tw_{1}+sw_{2})&=&J(tw_{1}+sw_{2})-\frac{1}{4}J'(tw_{1}+sw_{2})(tw_{1}+sw_{2})\\
&=&\biggl[\frac{1}{2}\widehat{M}(\|tw_{1}+sw_{2}\|^{2})-
\frac{1}{4}M(\|tw_{1}+sw_{2}\|^{2})\|tw_{1}+sw_{2}\|^{2}\biggl]\\
& +& \biggl[\displaystyle\int_{\Omega}\frac{1}{4}
f(tw_{1}+sw_{2})(tw_{1}+sw_{2})  - F(tw_{1}+sw_{2})\biggl] dx \\
&\leq &\biggl[\frac{1}{2}\widehat{M}(\|w_{1}+w_{2}\|^{2})-
\frac{1}{4}M(\|w_{1}+w_{2}\|^{2})\|w_{1}+w_{2}\|^{2}\biggl] 
\\
&+& \biggl[\displaystyle\int_{\Omega}\frac{1}{4}
f(w_{1}+w_{2})(w_{1}+w_{2})  - F(w_{1}+w_{2})\biggl] dx\leq\displaystyle\liminf_{n\rightarrow
+\infty} J(w_n) =c_{0}.
\end{eqnarray*}

To complete the proof of Theorem \ref{T2}, we claim that $w$ is a critical point for
functional $J$.  To see why, for each $\varphi \in H^{1}_{0}(\Omega)$, we introduce the functions 
$Q^{i}:\mathbb{R}^{3} \to \mathbb{R}$ given by
$$
\begin{array}{l}
Q^{1}(r,z,l)=M(\|w+r\varphi+zw_1+lw_2\|^{2})\|w+r\varphi+zw_1\|_1^{2} -\displaystyle \int_{\Omega_1}f(w_1+r\varphi_1+zw_1)(w_1+r\varphi_1+zw_1)dx.
\end{array}
$$
and
$$
\begin{array}{l}
Q^{2}(r,z,l)=M(\|w+r\varphi+zw_1+lw_2\|^{2})\|w+r\varphi+lw_2\|_2^{2} -\displaystyle \int_{\Omega_2}f(w_2+r\varphi_1+lw_2)(w_2+r\varphi_1+lw_2)dx.
\end{array}
$$
By a direct computation,
$$
\begin{array}{l}
\frac{\partial Q^{1}}{\partial z}(0,0,0)=2(M'(\|w\|^{2})\|w_1\|^{4}+M(\|w\|^{2})\|w_1\|^{2}) -\displaystyle \int_{\Omega_1}(f'(w_{1})w_1^{2} +f(w_1)w_1) dx.
\end{array}
$$
By inequality (\ref{relacaoderivadafuncaoM}) and (\ref{relacaoderivadafuncaof}), 
$$
-\frac{\partial Q^{1}}{\partial z}(0,0,0)>  2M'(\|w\|^{2})\|w_1\|^{2}\|w_2\|^{2}.
$$
Using a similar argument, it is possible to prove that
$$
-\frac{\partial Q^{2}}{\partial l}(0,0,0)>  2M'(\|w\|^{2})\|w_1\|^{2}\|w_2\|^{2} \quad \mbox{and} \quad \frac{\partial Q^{1}}{\partial l}(0,0,0)=\frac{\partial Q^{2}}{\partial z}(0,0,0)=M'(\|w\|^{2})\|w_2\|^{2}\|w_1\|^{2}.
$$
From this,
$$
\left|
\begin{array}{ll}
\frac{\partial Q^{1}}{\partial z}(0,0,0)   &  \frac{\partial Q^{2}}{\partial l}(0,0,0)\\
\frac{\partial Q^{1}}{\partial l}(0,0,0)   &  \frac{\partial Q^{2}}{\partial l}(0,0,0)
\end{array}
\right|=3(M'(\|w\|^{2}))^{2}\|w_2\|^{4}\|w_1\|^{4}>0.
$$
Therefore, applying the implicit function theorem, there are functions $z(r), l(r)$ of class $C^{1}$ defined on some interval $(-\delta, \delta), \delta>0$ such that $z(0)=l(0)=0$ and
$$
Q^{i}(r,z(r),l(r))=0, \,\,\, r \in (-\delta, \delta), i=1,2.
$$
This shows that for any $r \in (-\delta, \delta)$, 
$$
v(r)=w+r\varphi +z(r)w_1+l(r)w_2 \in {\mathcal M}.
$$
Then
$$
J(v(r)) \geq J(w), \,\,\,\, \forall r \in (-\delta, \delta),
$$
that is,
$$
J(w+r\varphi +z(r)w_1+l(r)w_2) \geq J(w), \,\,\,\, \forall r \in (-\delta, \delta).
$$
From this,
$$
\frac{J(w+r\varphi +z(r)w_1+l(r)w_2) - J(w) }{r} \geq 0, \,\,\,\, \forall r \in (0, \delta).
$$
Taking the limit of $r \to 0$, we get
$$
J'(w)(\varphi +z'(0)w_1+l'(0)w_2) \geq 0.
$$
Recalling that $J'(w)w_1=J'(w)w_2=0$, the above inequality loads to
$$
J'(w)\varphi \geq 0, \,\,\, \forall \varphi \in H^{1}_{0}(\Omega)
$$
and so,
$$
J'(w)\varphi= 0, \,\,\, \forall \varphi \in H^{1}_{0}(\Omega),
$$
showing that $w$ is a critical point for $J$. \qed

\section{An auxiliary Kirchhoff problem}

In this section, we work with an auxiliary problem adapting the ideas explored by del Pino \& Felmer in \cite{DelPinoFelmer} (see also \cite{Alves} and \cite{DingTanaka}).

We start recalling that the energy functional $ I_\lambda \colon E_\lambda \to \mathbb R $ associated with $ (P)_\lambda $ is  given by
$$
   I_\lambda (u) = \frac{1}{2} \widehat{M}(\| u \|_\lambda^{2})- \int_{ \mathbb R^3 } F(u)dx,
$$
where $ E_\lambda = \big( E, \| \cdot \|_\lambda \big) $ with
$$
   E = \left\{ u \in H^{1} ( \mathbb R^3 ) \, ; \, \int_{ \mathbb R^3 } a(x) |u|^{ 2 }dx < \infty \right\}
$$
and
$$
   \| u \|_\lambda = \left(\int_{\mathbb{R}^{3}}(|\nabla u|^{2}+(\lambda a(x)+1)|u|^{2})dx \right)^{\frac{1}{2}}.
$$
By $(a_0)$, the embedding below 
$$ 
E_\lambda \hookrightarrow H^{1}( \mathbb R^3 ) 
$$ 
is continuous for all $ \lambda \geq 0 $. Consequently, $ E_\lambda $ is compactly embedded in $ L_{ loc }^{s}( \mathbb R^3 ) $, for all $ 1 \leq s <6 $. A  direct computation gives 
that $ E_\lambda $ is a Hilbert space. Moreover, if $ {\cal O} \subset \mathbb R^3 $ is an open set, from the relation
\begin{equation} \label{modular relation 1}
  \int_{ \cal O } \left( \big| \nabla u \big|^{ 2} + \big( \lambda a(x) + 1 \big) | u |^{ 2 } \right)dx  \ge \int_{ \cal O } |u|^{ 2 }dx , \quad \forall u \in E_\lambda
\end{equation}
fixed $ \delta \in (0,1)$, there is $ \nu > 0 $, such that
\begin{equation} \label{modular relation 2}
   \|u\|_{ \lambda, \cal{O} }^{2} - \nu |u|_{2, \cal{O} }^{2} 
	\geq \delta \|u\|_{ \lambda,\cal{O} }^{2}, \, \forall u \in E_\lambda, \, \lambda \geq 0,
\end{equation}
where,
$$
\|u\|_{ \lambda, \cal{O} }= \left(\int_{ \cal O }( | \nabla u |^{ 2} + ( \lambda a(x) + 1 ) | u |^{ 2 })dx \right)^{\frac{1}{2}}
$$
and
$$
 |u|_{2, \cal{O} }= \left(\int_{ \cal O } | u |^{ 2 }dx \right)^{\frac{1}{2}}.
$$

We recall that given $ \epsilon > 0 $, the hypotheses $ (f_1) $ and $ (f_2) $ yield
\begin{equation} \label{f estimate}
   |f(s)| \le \epsilon | s | + C_\epsilon | s|^{5},  \,\,\, \mbox{and} \,\,\, s \in \mathbb R.
\end{equation}
Hence,
\begin{equation} \label{F estimate}
   |F(s)| \le \frac{\epsilon}{2} | s |^{2 } + \frac{C_\epsilon}{6} |s|^{6}, \, \forall   s  \in \mathbb R,
\end{equation}
where $ C_\epsilon $ depends on $ \epsilon $. Moreover, for  $\nu >0$ fixed in (\ref{modular relation 2}), the assumptions $ (f_1) $ and $ (f_4) $ imply that there is a unique $\xi>0$ verifying
\begin{equation} \label{NU}
\frac{f(\xi)}{\xi}=\nu
\end{equation}

Using the numbers $\xi$ and $\nu$, we set the function $ \tilde{f} \colon \mathbb R \to \mathbb R $ given by		 
$$
   \tilde{f}(s) =
   \begin{cases}
      \ \, f(s), \ s \le \xi \\
	    \nu \, s, \ s \ge \xi 
   \end{cases},
$$
which fulfills the inequality
\begin{equation} \label{til f estimate}
   \tilde{f}(s) \le \nu | s|, \,\,\, \forall  s \in \mathbb R.
\end{equation}
Thus
\begin{equation} \label{t til f estimate}
   \tilde{f}(s) s \le \nu | s |^{2}, \,\,\, \forall s \in \mathbb R
\end{equation}
and
\begin{equation} \label{til F estimate}
  \tilde{F}(s) \le \frac{\nu}{2} | s|^{ 2 }, \, \forall s \in \mathbb R,
\end{equation}
where $ \tilde{F}(s) = \int_0^t \tilde{f}(t) \, dt $.

Now, since $ \Omega=int( a^{ -1 } (\{0\})) $ is formed by $ k $ connected components $ \Omega_1, \ldots, \Omega_k $ with $ \text{dist} \big( \Omega_i, \Omega_j \big) > 0, \, i \ne j $, then for each $ j \in \{ 1, \ldots, k \} $, we are able to fix a smooth bounded domain $ \Omega'_j $ such that
\begin{equation} \label{omega}
   \overline{\Omega_j} \subset \Omega'_j \, \text{ and } \, \overline{\Omega'_i} \cap \overline{\Omega'_j} = \emptyset, \text{ for } i \ne j.
\end{equation}

From now on, we fix a non-empty subset $ \Upsilon \subset \left\{ 1, \ldots, k \right\}, $ 
$$
   \Omega_\Upsilon = \bigcup_{ j \in \Upsilon } \Omega_j, \, \Omega'_\Upsilon = \bigcup_{ j \in \Upsilon } \Omega'_j \,\, \mbox{and} \,\,\,
	 \chi_\Upsilon =
	 \begin{cases}
	    1, \text{ if } x \in \Omega'_\Upsilon \\
		  0, \text{ if } x \notin \Omega'_\Upsilon .
	 \end{cases}
$$
Using the above notations, we set the functions
$$
   g(x,s) = \chi_\Upsilon(x) f(s) + \big( 1-\chi_\Upsilon(x) \big) \tilde{f}(s), \, (x,s) \in \mathbb R^3 \times \mathbb R
$$
and
$$
   G(x,s) = \int_0^s g(x,t) \, dt, \, (x,s) \in \mathbb R^3 \times \mathbb R,
$$
and the auxiliary Kirchhoff problem
$$
\left\{
\begin{array}{l}
M(\|u\|^{2}_{\lambda})\biggl(- \Delta u + ( \lambda a(x) + 1 ) u\biggl)  = g(x,u), \text{ in } \mathbb R^3, \\
u \in E_\lambda .
\end{array}
\right.
\leqno{(A_\lambda)}
$$

The problem $ \big( A_\lambda \big) $ is strongly related to $ \big( P_\lambda \big) $, in the sense that, if $ u_\lambda $ is a solution for $ \big( A_\lambda \big) $  verifying
$$
    u_\lambda (x) \leq \xi, \, \forall x \in \mathbb R^N \setminus \Omega'_\Upsilon,
$$
then it is a solution for $ \big( P_\lambda \big) $.

In comparison to $ \big( P_\lambda \big) $, problem $ \big( A_\lambda \big) $ has the advantage that the energy functional associated with $ \big( A_\lambda \big) $, namely, $ \phi_\lambda \colon E_\lambda \to \mathbb R $ given by
$$
    \phi_\lambda(u) = \frac{1}{2} \widehat{M}(\|u\|^{2}_{\lambda})- \int_{ \mathbb R^3} G(x,u)dx,
$$
satisfies the $ (PS) $ condition, whereas $ I_\lambda $ does not necessarily satisfy this condition.

\vspace{0.5 cm}

\begin{proposition} \label{boundedness}
   All $ (PS)_d $ sequences for $ \phi_\lambda $ are bounded in $ E_\lambda $.
\end{proposition}

\noindent {\bf Proof.}
   Let $ (u_n) $ be a $ (PS)_d $ sequence for $ \phi_\lambda $. So, there is $ n_0 \in \mathbb N $ such that
$$
   \phi_\lambda (u_n) - \frac{1}{\theta} \phi_\lambda'(u_n) u_n \le d+1 + \| u_n \|_\lambda, \text { for } n \ge n_0.
$$
On the other hand, by (\ref{t til f estimate}) and (\ref{til F estimate}),
$$
   \tilde{F}(s) - \frac{1}{\theta} \tilde{f}(s)s \le \left( \frac{1}{2} - \frac{1}{\theta} \right) \nu | s |^{ 2 }, \, \forall x \in \mathbb R^3, s \in \mathbb R,
$$
which together with (\ref{modular relation 2}) gives
$$
   \phi_\lambda (u_n) - \frac{1}{\theta} \phi_\lambda'(u_n) u_n \ge \left( \frac{1}{2} - \frac{1}{\theta} \right) \delta \|u_n\|^{2}_{\lambda}, \, \forall n \in \mathbb N,
$$
from where it follows that $ (u_n) $ is bounded in $ E_\lambda $.
\qed

\begin{proposition} \label{Estimativa no infinito}
If $(u_n)$ is a $(PS)_d$ sequence for $\phi_{\lambda}$, then given $\epsilon>0$, there is $R>0$ such that
\begin{equation} \label{Estimativa}
   \limsup_n \int_{ \mathbb R^3 \setminus B_R (0) } (| \nabla u_n|^{2 } + ( \lambda a(x) + 1 ) | u_n |^{ 2 } )dx < \epsilon.
\end{equation}
Hence, once that $g$ has a subcritical growth, if $ u \in E_\lambda $ is the weak limit of $ (u_n) $, then
$$
\int_{\mathbb R^3}g(x,u_n)u_n\,dx \to \int_{\mathbb R^3} g(x,u)u \, dx \, \text{ and } \, \int_{\mathbb R^3} g(x,u_n)v \, dx \to \int_{\mathbb R^3} g(x,u)v \, dx, \, \forall v \in E_\lambda.
$$
\end{proposition}

\noindent {\bf Proof.}
   Let $ (u_n) $ be a $ (PS)_d $ sequence for $ \phi_\lambda $, $ R > 0 $ large such that $ \Omega'_\Upsilon \subset B_{ \frac{R}{2} }(0) $ and $ \eta_R \in C^\infty (\Bbb R^3) $ satisfying
$$
   \eta_R (x) =
	 \begin{cases}
	    0, \, x \in B_{ \frac{R}{2} }(0) \\
			1, \, x \in \Bbb R^3 \setminus B_R (0)
	 \end{cases},
$$
$ 0 \le \eta_R \le 1 $ and $ \big| \nabla \eta_R \big| \le \dfrac{C}{R} $, where $ C > 0 $ does not depend on $ R $. This way,
\begin{align*}
   \mbox{}  & m_{0}\|u_{n}\eta_{R}\|^{2}_{\lambda}\leq \int_{ \Bbb R^3 } M(\|u_{n}\|^{2}_{\lambda})(  | \nabla u_n \big|^{2 } 
	+ ( \lambda a(x) + 1 ) | u_n |^{ 2 } ) \eta_R dx\\
	   = & \phi_\lambda'(u_n) \left( u_n \eta_R \right) - \int_{ \Bbb R^3 }M(\|u_{n}\|^{2}_{\lambda}) u_n  \nabla u_n \nabla \eta_R dx+ \int_{ \Bbb R^3 \setminus \Omega'_\Upsilon } \tilde{f}(u_n) u_n \eta_Rdx.
\end{align*}
Denoting
$$
L =m_{0}\int_{ \Bbb R^3 } (  | \nabla u_n |^{ 2 } + ( \lambda a(x) +1 ) | u_n |^{2 } ) \eta_R dx,
$$
it follows from (\ref{t til f estimate}),
$$
   L \leq \phi_\lambda'(u_n) ( u_n \eta_R ) + \frac{C}{R} \int_{ \Bbb R^3 } M(\|u_{n}\|^{2}_{\lambda})| u_n | | \nabla u_n |dx + \nu \int_{ \Bbb R^3 } | u_n |^{2 } \eta_R dx.
$$
The H\"older's inequality in conjunction with the boundedness of $ (u_n) $ and $ ( | \nabla u_n | ) $ in $ L^{2}( \Bbb R^3 ) $, ensures that
$$
L \leq o_n(1) + \frac{C}{(1-\nu)R}.
$$
Therefore
$$
   \limsup_n \int_{ \Bbb R^3 \setminus B_R(0) } m_{0}(  | \nabla u_n |^{2 } + ( \lambda a(x) +1 ) | u_n |^{2} )dx \le \frac{C}{(1-\nu)R}.
$$
So, given $ \epsilon > 0 $,  choosing a $ R > 0 $ possibly still bigger, we have that $ \dfrac{C}{(1-\nu)R} < \epsilon $, which proves (\ref{Estimativa}). Now, we will show that
$$
\int_{\Bbb R^3}g(x,u_n)u_ndx \to \int_{\Bbb R^3}g(x,u)udx.
$$
Using the fact that $g(x,u)u \in L^{1}(\Bbb R^3)$ together with (\ref{Estimativa}) and Sobolev embeddings, given  $\epsilon >0$, we can choose $R>0$ such that
$$
\limsup_{n \to +\infty}\int_{\Bbb R^3 \setminus B_R(0)}|g(x,u_n)u_n|dx \leq \frac{\epsilon}{4} \quad \mbox{and} \quad \int_{\Bbb R^3 \setminus B_R(0)}|g(x,u)u|dx \leq \frac{\epsilon}{4}.
$$
On the other hand, since $g$ has a subcritical growth, the compact Sobolev embeddings load to
$$
\int_{B_{R}(0)}g(x,u_n)u_ndx \to \int_{B_{R}(0)}g(x,u)udx.
$$
Combining the above information, we conclude that
$$
\int_{\Bbb R^3}g(x,u_n)u_ndx \to \int_{\Bbb R^3}g(x,u)udx.
$$
The same type of argument works to prove that
$$
\int_{\Bbb R^3}g(x,u_n)v dx\to \int_{\Bbb R^3}g(x,u)vdx, \quad \forall v \in E_{\lambda}.
$$

\qed

The next result does not appear in \cite{DingTanaka}, however since we are
working with the Kirchhoff problem type, it is required here.

\begin{proposition} \label{GiovanyJoao}
  Let $(u_{n})$ be a $(PS)_{d}$ sequence for $\phi_{\lambda}$ such that
$u_{n}\rightharpoonup u$, then
$$
\lim_{n\to \infty}\int_{B_{R}}\left[ |\nabla u_{n}|^{2}+(\lambda a(x)+1)
u_{n}^{2}\right] dx=\int_{B_{R}}\left[ |\nabla u|^{2}+(\lambda a(x)+1)
u^{2}\right] dx,
$$
for all $R>0$.
\end{proposition}
		
\noindent {\bf Proof.} \, We can assume that $\|u_{n}\|_{\lambda}\rightarrow t_{0}$, thus, we have $\|u\|_{\lambda}\leq t_{0}$. Let $\eta_{\rho} \in C^{\infty}(\R^3)$ such that
\begin{eqnarray*}
\eta_{\rho}(x) =\ \  \left\{
             \begin{array}{l}
              1 \quad se \quad x \in B_{\rho}(0)
        \\
        \\
        0 \quad se \quad x \not\in B_{2\rho}(0).
        \\
        \\
             \end{array}
           \right.
\end{eqnarray*}
with $0 \leq \eta_{\rho}(x) \leq 1$. Let,
$$
P_{n}(x)=M(\|u_{n}\|_{\lambda}^{2})\left[ |\nabla u_{n}-\nabla u|^{2}+(\lambda a(x)+1)(u_{n}-u)^{2}\right].
$$
For each $R>0$ fixed, choosing $\rho>R$ we obtain
$$
\int_{B_{R}}P_{n}=\int_{B_{R}}P_{n}\eta_{\rho}\leq
M(\|u_{n}\|_{\lambda}^{2})\int_{\R^{3}}\left[ |\nabla u_{n}-\nabla
u|^{2}+(\lambda a(x)+1)(u_{n}-u)^{2}\right]\eta_{\rho}.
$$
By expanding the inner product in $\R^{3}$,
\begin{eqnarray*}
\int_{B_{R}}P_{n}&\leq& M(\|u_{n}\|_{\lambda}^{2})\int_{\R^{3}}\left[ |\nabla u_{n}|^{2}+(\lambda a(x)+1)(u_{n})^{2}\right]\eta_{\rho}\\
&-&2M(\|u_{n}\|_{\lambda}^{2})\int_{\R^{3}}\left[ \nabla u_{n}\nabla u+(\lambda a(x)+1)u_{n}u\right]\eta_{\rho}\\
&+&M(\|u_{n}\|_{\lambda}^{2})\int_{\R^{3}}\left[ |\nabla u|^{2}+(\lambda a(x)+1)u^{2}\right]\eta_{\rho}.
\end{eqnarray*}
Setting
$$
I_{n,\rho}^{1}=M(\|u_{n}\|_{\lambda}^{2})\int_{\R^{3}}\left[ |\nabla u_{n}|^{2}+(\lambda a(x)+1)(u_{n})^{2}\right]\eta_{\rho}
-\int_{\R^{3}}g( x, u_{n})u_{n}\eta_{\rho},
$$
$$
I_{n,\rho}^{2}=M(\|u_{n}\|_{\lambda}^{2})\int_{\R^{3}}\left[ \nabla u_{n}\nabla u+(\lambda a(x)+1)u_{n}u\right]\eta_{\rho}
-\int_{\R^{3}}g( x, u_{n})u\eta_{\rho},
$$
$$
I_{n,\rho}^{3}=-M(\|u_{n}\|_{\lambda}^{2})\int_{\R^{3}}\left[ \nabla u_{n}\nabla u+(\lambda a(x)+1)u_{n}u\right]\eta_{\rho}
+M(\|u_{n}\|_{\lambda}^{2})\int_{\R^{3}}\left[ |\nabla u|^{2}+(\lambda a(x)+1)u^{2}\right]\eta_{\rho}
$$
and
$$
I_{n,\rho}^{4}=\int_{\R^{3}}g( x, u_{n})u_{n}\eta_{\rho}-\int_{\R^{3}}g( x, u_{n})u\eta_{\rho},
$$
we find the estimate below
\begin{equation}\label{11}
0\leq\int_{B_{R}}P_{n}\leq |I_{n,\rho}^{1}|+|I_{n,\rho}^{2}|+|I_{n,\rho}^{3}|+|I_{n,\rho}^{4}|.
\end{equation}
Observe that
$$
I_{n,\rho}^{1}=\phi_{\lambda}'(u_{n})(u_{n}\eta_{\rho})-M(\|u_{n}\|_{\lambda}^{2})\int_{\R^{3}} u_{n}\nabla u_{n}\nabla\eta_{\rho}.
$$
Recalling that $(u_{n}\eta_{\rho})$ is bounded in $H_{\varepsilon}$, we have
$\phi_{\lambda}'(u_{n})(u_{n}\eta_{\rho})=o_{n}(1)$. Moreover, from a
straightforward computation
$$
\lim_{\rho\to \infty}\left[ \limsup_{n\to \infty}
\left|M(\|u_{n}\|_{\lambda}^{2})\int_{\R^{3}} u_{n}\nabla
u_{n}\nabla\eta_{\rho}\right|\right]=0.
$$
Then,
\begin{equation}\label{12}
\lim_{\rho\to \infty}\left[ \limsup_{n\to \infty} |I_{n,\rho}^{1}|\right]=0.
\end{equation}
We also see that
$$
I_{n,\rho}^{2}=\phi_{\lambda}'(u_{n})(u\eta_{\rho})-M(\|u_{n}\|_{\lambda}^{2})\int_{\R^{3}} u\nabla u_{n}\nabla\eta_{\rho}.
$$
By arguing in the same way as in the previous case,
$$
\phi_{\lambda}'(u_{n})(u\eta_{\rho})=o_{n}(1)
$$
and
$$
\lim_{\rho\to \infty}\left[ \limsup_{n\to \infty} \left|M(\|u_{n}\|_{\lambda}^{2})\int_{\R^{3}} u\nabla u_{n}\nabla\eta_{\rho}\right|\right]=0.
$$
Therefore,
\begin{equation}\label{13}
\lim_{\rho\to \infty}\left[ \limsup_{n\to \infty} |I_{n,\rho}^{2}|\right]=0.
\end{equation}
On the other hand, from the weak convergence
\begin{equation}\label{14}
\lim_{n\to \infty} |I_{n,\rho}^{3}|=0, \ \forall \ \rho>R.
\end{equation}
Finally, from
$$
u_{n}\rightarrow u, \ \mbox{in} \ L^{s}_{loc}(\R^{3}), 1\leq s<6,
$$
we conclude that
\begin{equation}\label{15}
\lim_{n\to \infty} |I_{n,\rho}^{4}|=0, \ \forall \ \rho>R.
\end{equation}
From $(\ref{11})$, $(\ref{12})$, $(\ref{13})$, $(\ref{14})$ and $(\ref{15})$, we obtain
$$
0\leq \limsup_{n\to \infty} \int_{B_{R}}P_{n}\leq 0,
$$
or equivalently 
$$
\displaystyle\lim_{n\to \infty} \displaystyle\int_{B_{R}}P_{n}=0.
$$ 
Therefore, 
$$
\lim_{n\to \infty}\int_{B_{R}}\left[ |\nabla u_{n}|^{2}+(\lambda a(x)+1) u_{n}^{2}\right]=\int_{B_{R}}\left[ |\nabla u|^{2}+(\lambda a(x)+1) u^{2}\right].
$$

\qed

\begin{proposition} \label{(PS) condition}
   $ \phi_\lambda $ verifies the $ (PS) $ condition.
\end{proposition}
		
\noindent {\bf Proof}
Let $ (u_n) $ be a $ (PS)_d $ sequence for $ \phi_\lambda $ and $ u \in E_\lambda $ such that $u_n \rightharpoonup u$ in $E_{\lambda}$. Thereby, by Proposition \ref{Estimativa no infinito},
$$
   \int_{ \mathbb R^3 } g(x,u_n)u_ndx \to \int_{ \mathbb R^3 } g(x,u)u dx\, \text{ and } \, \int_{ \mathbb R^3 } g(x,u_n)v dx\to \int_{ \mathbb R^3 } g(x,u)vdx, \, \forall v \in E_\lambda.
$$
Moreover, the weak limit also gives
$$
\int_{ \mathbb R^3 }  \nabla u \nabla ( u_n-u )dx  \to 0
$$
and
$$
\int_{ \mathbb R^3 }  ( \lambda a(x) + 1 ) u ( u_n-u )dx  \to 0.
$$
Gathering $\phi_\lambda'(u_n)u_n=o_n(1)$, $\phi_\lambda'(u_n)u=o_n(1)$, $(M_{1})$ and the above limits, we derive that 
$$
 \|u_n-u\|^{2}_\lambda \to 0,
$$
finishing the proof. \qed


\section{The $ (PS)_\infty $ condition}


A sequence $ (u_n) \subset  H^{ 1} ( \mathbb R^3 ) $ is called a $ (PS)_\infty $ \emph{sequence for the family} $ \left( \phi_\lambda \right)_{\lambda \ge 1} $, if there is a sequence $ ( \lambda_n ) \subset  [1, \infty) $ with $ \lambda_n \to \infty $, as $ n \to \infty $, verifying
$$
   \phi_{ \lambda_n }(u_n) \to c \text{ and } \left\| \phi'_{ \lambda_n }(u_n) \right\|_{E^{*}_{\lambda_n}} \to 0, \text{ as } n \to \infty,
$$
for some $c \in \mathbb{R}$.
\begin{proposition} \label{(PS) infty condition}
   Let $ (u_n) \subset H^{ 1, } ( \mathbb R^3 ) $ be a $ (PS)_\infty $ sequence for $ \left( \phi_\lambda \right)_{\lambda \ge 1} $. Then, up to a subsequence, there exists $ u \in H^{ 1 } ( \mathbb R^3 ) $ such that $ u_n \rightharpoonup u $ in $ H^{ 1}( \mathbb R^3 ) $. Furthermore,
\begin{enumerate}
  \item[(i)] $ u_n \to u $ in $ H^{ 1} ( \mathbb R^3 ) $;
	\item[(ii)] $ u = 0 $ in $ \mathbb R^3 \setminus \Omega_\Upsilon $,  $ u_{|_{\Omega_j}} \geq 0$ for all $j \in \Upsilon $, and $u$ is a solution for
   $$
   	    \begin{cases}
		     M\biggl(\displaystyle\int_{\Omega_\Upsilon}(\nabla u|^{2}+u^{2})dx\biggl)\biggl(- \Delta u + u \biggl) = f(u), \text{ in } \Omega_\Upsilon, \\
				   u \in H^{1}_0 ( \Omega_\Upsilon );
			\end{cases}
   \eqno{(P)_{\infty,\Upsilon} }
	$$
	 \item[(iii)] $ \displaystyle \lambda_n\int_{\mathbb R^3}  a(x) | u_n |^{ 2 }dx \to 0 $;
	 \item[(iv)] $ \|u_n-u\|^{2}_{\lambda,\Omega'_\Upsilon}\to 0, \text{ for } j \in \Upsilon $;
   \item[(v)] $  \|u_n\|^{2}_{\lambda,\mathbb{R}^{3} \setminus \Omega'_\Upsilon} \to 0 $;
	 \item[(vi)] $ \phi_{ \lambda_n } (u_n) \to \displaystyle \frac{1}{2}\widehat{M}\biggl(\int_{ \Omega_\Upsilon }(| \nabla u |^{ 2 } +  
	| u |^{2})dx\biggl) 	- \int_{ \Omega_\Upsilon } F(u)dx $.
\end{enumerate}
\end{proposition}

\noindent {\bf Proof.}
   Using the Proposition \ref{boundedness}, we know that $( \| u_n \|_{ \lambda_n } ) $ is bounded in $ \mathbb R $ and $ (u_n) $ is bounded in $ H^{1}( \mathbb R^3 ) $. So, up to a subsequence, there exists $ u \in H^{1}(\mathbb R^3) $ such that
$$
   u_n \rightharpoonup u  \text{ in } H^{1} ( \mathbb R^3) \, \text{ and } \, u_n(x) \to u(x) \text{ for a.e. } x \in \mathbb R^3.
$$
Now, for each $ m \in \mathbb N $, we define $ C_m = \left\{ x \in \mathbb R^3 \, ; \, a(x) \ge \dfrac{1}{m} \right\} $. Without loss of generality, we can assume $ \lambda_n < 2 ( \lambda_n-1 ), \, \forall n \in \mathbb N $. Thus
	       $$
		       \int_{ C_m } | u_n |^{ 2 }dx \le \frac{2m}{\lambda_n} \int_{ C_m } \big( \lambda_n a(x)+1) | u_n |^{ 2 }dx \leq \frac{C}{\lambda_n}.
	       $$
By Fatou's lemma, 
	       $$
			     \int_{ C_m } | u |^{ 2 }dx = 0,
	       $$ 
implying that $ u = 0 $ in $ C_m $, and so, $ u = 0 $ in $ \mathbb R^3 \setminus \overline{\Omega} $. From this, we are able to prove $(i)-(vi)$.
	
\begin{enumerate}
   \item[$(i)$] Since $ u = 0 $  in $ \mathbb R^3 \setminus \overline{\Omega} $, repeating the argument explored in  Proposition \ref{(PS) condition}, we get
	       $$
            \int_{ \mathbb R^3 } (|\nabla u_n - \nabla u|^{2}  + ( \lambda_n a(x) + 1) |u_n-u|^{2})dx \to 0,
         $$
  which implies $ u_n \to u $ in $ H^{1}(\mathbb R^3) $.
	 \item[$(ii)$] Since $ u \in H^{ 1}(\mathbb R^3) $ and $ u = 0 $  in $ \mathbb R^3 \setminus \overline{\Omega} $, 
	we have $ u \in H^{1}_0( \Omega ) $ or, equivalently, $ u_{ |_{\Omega_j} } \in H^{1}_0( \Omega_j) $, for $ j = 1, \ldots, k $. 
	Moreover, the limit $u_n \to u$ in $H^{1}(\mathbb R^3)$ combined with $\phi'_{\lambda_n}(u_n)\varphi \to 0$  
	for $\varphi \in C^{\infty}_0 ( \Omega_\Upsilon )$ implies that
         \begin{equation} \label{u is solution}	
	         M\biggl(\int_{\Omega_\Upsilon}( |\nabla u|^{2}  +  u^{2} )dx\biggl) \int_{\Omega_\Upsilon}( \nabla u  \nabla \varphi 
					+  u \varphi )dx  - \int_{\Omega_\Upsilon} f(u) \varphi dx = 0,
         \end{equation}
showing that $ u_{ |_{\Omega_{\Upsilon}} }  $ is a solution for the nonlocal problem
         $$
				 \begin{cases}
		        M\biggl(\displaystyle \int_{\Omega_\Upsilon}( |\nabla u|^{2}  +  u^{2} )dx\biggl)(- \Delta u + u) = f(u), \text{ in } \Omega_\Upsilon, \\
		        u \in H^{1}_0 ( \Omega_\Upsilon ).
		       \end{cases}
					\eqno{(P)_{\infty,\Upsilon} }
$$
	
 On the other hand, if $ j \notin \Upsilon $, we must have
				 $$
				   M\biggl(\int_{\Omega_\Upsilon}( |\nabla u|^{2}  +  u^{2} )dx\biggl)\int_{ \Omega_j }(|\nabla u|^{2} + |u|^{2}) dx
					- \int_{ \Omega_j } \tilde{f}(u)u dx= 0.
				 $$
				The above equality combined with (\ref{t til f estimate}) and (\ref{modular relation 2}) gives
			 	 $$
				    0 \ge \|u\|^{2}_{ \lambda, \Omega_j }- \nu \|u\|^{2}_{2,\Omega_j } \geq \delta \|u\|^{2}_{ \lambda, \Omega_j }(u) \geq 0,
				 $$
				 from where it follows $ u_{|_{\Omega_j}} = 0 $ for $ j \notin \Upsilon $.  This proves $ u = 0 $ outside $ \Omega_\Upsilon $ and $ u \geq 0 $ in $ \Bbb R^3$.
				
	\item[$(iii)$]  From (i),
	       $$
		        \lambda_n \int_{ \mathbb R^3 }  a(x) | u_n |^{2 } dx= \int_{ \mathbb R^2 } \lambda_n a(x) | u_n-u |^{ 2 }dx \leq \|u_n-u\|^{2}_{\lambda_n},
	       $$
				loading to
				$$
				\lambda_n \int_{ \mathbb R^3 }  a(x) | u_n |^{2 } dx \to 0.
				$$
	\item[$(iv)$] Let $ j \in \Upsilon $. From (i),
	       $$
				    |u_n-u|^{2}_{2, \Omega'_j }, |\nabla u_n - \nabla u|^{2}_{2, \Omega'_j } \to 0.
				 $$
				 Then,
				 $$
				 \int_{ \Omega'_\Upsilon } ( | \nabla u_n |^{ 2 } - | \nabla u |^{ 2 } )dx \to 0 \quad \mbox{and} \quad \int_{ \Omega'_\Upsilon } ( | u_n |^{ 2 } - | u |^{ 2 } )dx \to 0.
				 $$
				 From (iii),
				 $$
				    \int_{ \Omega'_\Upsilon } \lambda_n a(x) | u_n |^{2}  dx \to 0.
			   $$
				 This way
				 $$
				  \|u_n\|^{2}_{ \lambda_n, \Omega'_\Upsilon } \to \int_{ \Omega_\Upsilon } ( | \nabla u |^{ 2 } + | u |^{ 2} )dx.
				 $$
	 \item[$(v)$] By (i),  $ \|u_n-u\|^{2}_{ \lambda_n } \to 0 $, and so,
	              $$
								  \|u_n\|^{2}_{ \lambda_n, \mathbb R^3 \setminus \Omega_\Upsilon } \to 0.
								$$
	 \item[$(vi)$] From $(i)-(v)$,
				 $$
					  \widehat{M}\biggl(\|u_n\|^{2}_{\lambda_n}\biggl)
						\to \widehat{M}\biggl(\int_{ \Omega_{\Upsilon} }(| \nabla u|^{2} + | u |^{2})dx\biggl)
				 $$
				and
				 $$
					  \int_{ \mathbb R^3 } G(x,u_n)dx \to \int_{ \Omega_\Upsilon } F(u)dx.
				 $$
				 Therefore
				 $$
					\phi_{ \lambda_n } (u_n) \to \frac{1}{2}\widehat{M}\biggl(\int_{ \Omega_\Upsilon } ( | \nabla u |^{2} + | u |^{2} )dx\biggl)
					- \int_{ \Omega_\Upsilon } F(u)dx.
				 $$
\end{enumerate}	
\qed


\section{The boundedness of the $ \big( A_\lambda \big) $ solutions}


In this section, we study the boundedness outside $ \Omega'_\Upsilon $ for some solutions of $ \big( A_\lambda \big ) $. To this end, we adapt the arguments found in \cite{Alves} and \cite{Gon} for our new setting.

\begin{proposition} \label{P:boundedness of the solutions}
    Let $ \big( u_\lambda \big) $ be a family of solutions for $ \big( A_\lambda \big) $ such that $ u_\lambda \to 0 $ in $ H^{1}( \mathbb R^3 \setminus \Omega_\Upsilon ) $, as $ \lambda \to \infty $. Then, there exists $ \lambda^* > 0 $ with the following property:
$$
   \left| u_\lambda \right|_{ \infty, \mathbb R^3 \setminus \Omega'_\Upsilon } \le \xi, \, \forall \lambda \ge \lambda^*.
$$
Hence, $u_{\lambda}$ is a solution for $(P_\lambda)$ for $\lambda \geq \lambda^*$.
\end{proposition}

\vspace{.3cm}
\noindent {\bf Proof.}
Since $\partial \Omega'_\Upsilon$ is a compact set, fixed a neighborhood $\mathcal B$ of $\partial \Omega'_\Upsilon$ such that
$$
\mathcal B \subset \mathbb{R}^{3} \setminus \Omega_\Upsilon,
$$
the interation Moser technique implies that there is $C>0$, which is independent of $\lambda$, such that
$$
|u_\lambda|_{L^{\infty}(\partial \Omega'_\Upsilon) }
 \leq C |u_{\lambda}|_{L^{2^{*}}({\mathcal B})}
$$
Since $u_\lambda \to 0$ in $H^{1}(\mathbb{R}^{3} \setminus \Omega_\Upsilon)$, we have that
$$
|u_{\lambda}|_{L^{2^{*}}({\mathcal B})} \to 0 \,\,\, \mbox{as} \,\,\, \lambda \to \infty.
$$
Hence, there is $\lambda^{*}>0$ such that
$$
|u_{\lambda}|_{L^{2^{*}}({\mathcal B})} < \frac{\xi}{2C} \,\,\, \forall \lambda \geq \lambda^{*},
$$
and so,
$$
|u_\lambda|_{L^{\infty}(\partial \Omega'_\Upsilon) }< \xi \,\,\, \forall \lambda \geq \lambda^{*}.
$$
Next, for $ \lambda \ge \lambda^* $, we set  $ \widetilde{u}_\lambda \colon \mathbb R^3 \setminus \Omega'_\Upsilon \to \mathbb R $ given by
$$
   \widetilde{u}_\lambda(x) = ( u_\lambda-\xi)^+ (x).
$$
Thereby,  $ \widetilde{u}_\lambda \in H^{1}_0(\mathbb R^3 \setminus \Omega'_\Upsilon) $. Our goal is showing that  $\widetilde{u}_\lambda = 0 $ in $ \mathbb R^3 \setminus \Omega'_\Upsilon $, because this will imply that 
$$
   \left| u_\lambda \right|_{ \infty, \mathbb R^3 \setminus \Omega'_\Upsilon } \leq \xi.
$$
In fact, extending $ \widetilde{u}_\lambda = 0 $ in $ \Omega'_\Upsilon $ and taking $ \widetilde{u}_\lambda $ as a test function 
and using $(M_{1})$, we obtain
$$
   m_{0}\biggl(\int_{ \mathbb R^3 \setminus \Omega'_\Upsilon } \nabla u_\lambda  \nabla \widetilde{u}_\lambda dx 
	+ \int_{ \mathbb R^3 \setminus \Omega'_\Upsilon } \! \! \! \! ( \lambda a(x) + 1)u_\lambda \widetilde{u}_\lambda dx  \biggl)\leq \int_{ \mathbb R^N \setminus \Omega'_\Upsilon } g \left( x, u_\lambda \right) \widetilde{u}_\lambda dx.
$$
Since
\begin{gather*}
   \int_{ \mathbb R^3 \setminus \Omega'_\Upsilon } \nabla u_\lambda \nabla \widetilde{u}_\lambda dx = \int_{ \mathbb R^3 \setminus \Omega'_\Upsilon }| \nabla \widetilde{u}_\lambda |^{2} dx, \\
	\int_{ \mathbb R^3 \setminus \Omega'_\Upsilon } ( \lambda a(x) +1) u_\lambda \widetilde{u}_\lambda dx= \int_{ \left( \mathbb R^3 \setminus \Omega'_\Upsilon \right)_+ } ( \lambda a(x) + 1) \left( \widetilde{u}_\lambda+\xi \right) \widetilde{u}_\lambda dx
\end{gather*}
and
$$
  \int_{ \mathbb R^3 \setminus \Omega'_\Upsilon } g \left( x, u_\lambda \right) \widetilde{u}_\lambda dx = \int_{ \left( \mathbb R^3 \setminus \Omega'_\Upsilon \right)_+ } \frac{g \left( x, u_\lambda \right)}{u_\lambda} \left( \widetilde{u}_\lambda+\xi \right) \widetilde{u}_\lambda dx,
$$
where
$$
   \left( \mathbb R^3 \setminus \Omega'_\Upsilon \right)_+ = \left\{ x \in \mathbb R^3 \setminus \Omega'_\Upsilon \, ; \, u_\lambda(x) > \xi \right\},
$$
we derive
$$
   m_{0}\biggl(\int_{ \mathbb R^3 \setminus \Omega'_\Upsilon } | \nabla \widetilde{u}_\lambda|^{2 }dx 
	+ \int_{ \left( \mathbb R^3 \setminus \Omega'_\Upsilon \right)_+ } \! \! \! \! ( (\lambda a(x) + 1 )\biggl)-\frac{g \left( x, u_\lambda \right)}{u_\lambda}) \left( \widetilde{u}_\lambda+\xi \right)  \widetilde{u}_\lambda dx \leq 0,
$$
Now, by \eqref{til f estimate},
$$
 ( \lambda a(x) + 1) - \frac{g ( x, u_\lambda )}{u_\lambda} > \nu  - \frac{\tilde{f} \left( x, u_\lambda \right)}{u_\lambda} \ge 0 \quad \mbox{in} \quad  \left( \mathbb R^3 \setminus \Omega'_\Upsilon \right)_+ .
$$
This form, $ \widetilde{u}_\lambda = 0 $ in $ \left( \mathbb R^3 \setminus \Omega'_\Upsilon \right)_+ $. Obviously, $ \widetilde{u}_\lambda = 0  $ at the points where $ u_\lambda \leq \xi $, consequently, $ \widetilde{u}_\lambda = 0 $ in $ \mathbb R^N \setminus \Omega'_\Upsilon $. 		
\qed


\section{A special minimax value for $ \phi_\lambda $}


For fixed non-empty subset $ \Upsilon \subset \left\{ 1, \ldots, k \right\} $, consider
$$
   I_\Upsilon(u) = \frac{1}{2}\widehat{M}\biggl(\int_{\Omega_\Upsilon} (| \nabla u |^{2} + | u |^{2} ) dx\biggl)- \int_{ \Omega_\Upsilon } F(u)dx, \ u \in H^{1}_0( \Omega_\Upsilon), 
$$
the energy functional associated to $ (P)_{\infty,\Upsilon} $,  and $\phi_{ \lambda,\Upsilon }: H^{1}(\Omega'_\Upsilon) \to \mathbb{R} $ given by 
$$
   \phi_{ \lambda,\Upsilon }(u) = \frac{1}{2}\widehat{M}\biggl(\int_{ \Omega'_\Upsilon } (| \nabla u |^{2} 
	+( \lambda a(x) + 1) | u |^{2} )dx \biggl) - \int_{ \Omega'_\Upsilon } F(u)dx, 
$$
the energy functional associated to the Kirchhoff problem 
$$
   \begin{cases}
	    M\biggl(\displaystyle \int_{ \Omega'_\Upsilon } (| \nabla u |^{2} 
	+( \lambda a(x) + 1) | u |^{2} )dx \biggl)(- \Delta u + ( \lambda a(x) + 1) u ) = f(u), \text{ in } \Omega'_\Upsilon, \\
			\frac{\partial u}{\partial \eta} = 0, \text{ on } \partial \Omega'_\Upsilon.
	 \end{cases}
$$

In the following, we denote by $c_\Upsilon$ the number given by
$$
c_\Upsilon=\inf_{u \in \mathcal M_{\Upsilon}}I_\Upsilon(u)
$$
where
$$
\mathcal M_{\Upsilon}=\{u\in \mathcal N_{\Upsilon}: I_\Upsilon'(u)u_j=0 \mbox { and }
u_{j}\neq 0 \,\,\, \forall j \in \Upsilon \}
$$
with $u_{j}=u_{|_{ \Omega_j}}$ and
$$
\mathcal{N}_{\Upsilon} = \{ u\in
H^1_0(\Omega_\Upsilon)\setminus\{0\}\, :\, I_\Upsilon'(u)u=0\}.
$$
Of a similar way,  we denote by $c_{\lambda,\Upsilon}$ the number given by
$$
c_{\lambda,\Upsilon}=\inf_{u \in \mathcal M'_{\Upsilon}}\phi_{\lambda,\Upsilon}(u)
$$
where
$$
\mathcal M'_{\Upsilon}=\{u\in \mathcal N'_{\Upsilon}: \phi_{\lambda,\Upsilon}'(u)\tilde{u}_j=0 \mbox { and }
\tilde{u}_{j}\neq 0 \,\,\, \forall j \in \Upsilon \}
$$
with 
$$
\tilde{u}_j(x)=
\left\{
\begin{array}{l}
u(x), \quad x \in \Omega'_j\\
0, \quad x \in \Omega'_\Upsilon \setminus \Omega'_j
\end{array}
\right.
$$
and 
$$
\mathcal{N}'_{\Upsilon} = \{ u\in H^1(\Omega'_\Upsilon)\setminus\{0\}\, :\,\phi_{\lambda,\Upsilon}'(u)u=0\}.
$$
Repeating the same approach used in Section 2,  we ensure that there exist $ w_\Upsilon  \in H^{1}_0( \Omega_\Upsilon  ) $ and $ w_{ \lambda,\Upsilon  } \in H^{1}( \Omega'_\Upsilon ) $ such that
$$
   I_\Upsilon( w_\Upsilon) = c_\Upsilon  \, \text{ and } \, I'_\Upsilon ( w_\Upsilon) = 0
$$
and
$$
	\phi_{ \lambda,\Upsilon}( w_{ \lambda,\Upsilon }) = c_{ \lambda,\Upsilon}	\, \text{ and } \, \phi'_{ \lambda,\Upsilon }( w_{ \lambda,\Upsilon }) = 0.
$$

By a direct computation, it is possible to show that there is $\tau>0$ such that if $u \in \mathcal M_{\Upsilon}$, then
\begin{equation} \label{tau0}
\|u_j\|_j > \tau, \,\,\, \forall j \in \Upsilon,
\end{equation}
where, $\|\,\,\,\,\|_j$ denotes the norm on $H^{1}_0(\Omega_j)$ given by
$$
||u||_j=\left(\int_{\Omega_j}(|\nabla u|^{2}+|u|^{2})dx\right)^{\frac{1}{2}}.
$$
In particular, since $w_\Upsilon \in \mathcal M_{\Upsilon}$, we also have
\begin{equation} \label{tau1}
\|w_{\Upsilon,j}\|_j > \tau \,\,\, \forall j \in \Upsilon,
\end{equation}
where $w_{\Upsilon,j}=w_{\Upsilon}|_{\Omega_j}$ for all $j \in \Upsilon$. Moreover, reviewing the proof of Theorem \ref{T2}, it is possible to see that 
\begin{equation} \label{C01}
I_\Upsilon(w_\Upsilon)=c_\Upsilon=\max\{I_\Upsilon(t_1w_1+....+t_lw_l)\,:\,(t_1,....,t_l) \in [0,+\infty)^{l}\}
\end{equation}
and
\begin{equation} \label{C02}
I_\Upsilon(w_\Upsilon)=I_\Upsilon(w_1+....+w_l)>I_\Upsilon(t_1w_1+....+t_lw_l), \quad \forall  (t_1,....,t_l) \in [0,+\infty)^{l} \setminus \{(1,....,1)\}.
\end{equation}

\vspace{.3cm}
\begin{lemma}
   There holds that
\begin{enumerate}
   \item[(i)] $ 0 < c_{ \lambda,\Upsilon  } \le c_\Upsilon , \, \forall \lambda \geq 0 $;
   \item[(ii)] $ c_{ \lambda,\Upsilon  } \to c_\Upsilon , \text{ as } \lambda \to \infty$.
\end{enumerate}
\end{lemma}

\noindent {\bf Proof}
\begin{enumerate}
   \item[(i)]  Using the inclusion $ H^{1}_0(\Omega_\Upsilon) \subset H^{1}(\Omega'_\Upsilon) $, it  easy to observe that
				 $$
				   c_{ \lambda,\Upsilon  } \leq  c_\Upsilon .
				 $$
	 \item[(ii)]  Let $ \left( \lambda_n \right) $ be such a sequence with $\lambda_n \to +\infty$ and consider an arbitrary subsequence of $ \left( c_{ \lambda_n,\Upsilon  } \right) $ (not relabelled) . Let $ w_n \in H^{1}( \Omega'_j) $ with 
		     $$
					  \phi_{ \lambda_n,\Upsilon  }(w_n) = c_{ \lambda_n,\Upsilon }	\, \text{ and } \, \phi'_{ \lambda_n,\Upsilon } ( w_n) = 0.
				 $$
By the previous item, $ \big( c_{ \lambda_n,\Upsilon  } \big) $ is bounded. Then, there exists $( w_{ n_k }) $ subsequence of $( w_n) $ such that $ (\phi_{ \lambda_{ n_k },\Upsilon }( w_{ n_k }) )$ converges and  $\phi'_{ \lambda_{ n_k },\Upsilon } ( w_{ n_k }) = 0 $. Now, repeating the arguments explored in the proof of  Proposition \ref{(PS) infty condition}, there is $ w \in H^{1}_0(\Omega_\Upsilon ) \setminus \{0\} \subset H^{1}(\Omega'_\Upsilon ) $ such that
$$
w_j =w|_{\Omega_j} \not=0  \,\,\, \forall j \in \Upsilon
$$	
and 				
				
				$$
				    w_{ n_k } \to w \text{ in } H^{1}( \Omega'_\Upsilon ), \text{ as } k \to \infty.
				 $$
				 Furthermore, we also can prove that
				 $$
				    c_{ \lambda_{ n_k },\Upsilon } = \phi_{ \lambda_{ n_k },\Upsilon  }( w_{ n_k }) \to I_\Upsilon (w)
				 $$
				and
				 $$
				    0 = \phi'_{ \lambda_{ n_k },\Upsilon  }( w_{ n_k }) \to I'_\Upsilon (w).
				 $$
				Then,  $ w \in \mathcal M_{\Upsilon}$, and  by definition of $c_{\Upsilon}$,
				 $$
				    \lim_k c_{ \lambda_{ n_k },\Upsilon  } \geq c_\Upsilon .
				 $$
				 The last inequality together with item (i) implies
				 $$
				    c_{ \lambda_{ n_k },\Upsilon  } \to c_\Upsilon , \text{ as } k \to \infty.
				 $$
				 This establishes the asserted result.
\end{enumerate}
\qed
\vspace{.5cm}

In the sequel,  we fix  $ R > 1 $ verifying
\begin{equation} \label{R}
   0< I'_j \left( \frac{1}{R}w_j+\sum_{k=1, k\not=j}^{l}t_k R w_k \right)\left( \frac{1}{R} w_j \right) \,\,\, \mbox{and} \,\,\, I'_\Upsilon\left(R w_j+\sum_{k=1, k\not=j}^{l}t_k R w_k\right)(R w_j)<0,  
\end{equation}
for $ j \in \Upsilon$ and $\forall t_k \in [1/R^2,1] $ with $k \not=j$.

In the sequel, to simplify the notation, we rename the components $ \Omega_j $ of $ \Omega $ in way such that $ \Upsilon = \{ 1, 2, \ldots, l \} $  for some $ 1 \le l \le k $. Then, we define:
\begin{gather*}
   \gamma_0 ( \textbf{t} )(x) = \sum_{j=1}^l t_j R w_j(x) \in H^{1}_0(\Omega_\Upsilon), \, \forall \textbf{t}=( t_1, \ldots, t_l )\in [1/R^2,1]^l, \\
   \Gamma_\ast = \Big\{ \gamma \in C \big( [1/R^2,1]^l, E_\lambda \setminus \{ 0 \} \big) \, ; \, \gamma(\textbf{t})|_{\Omega'_j} \not=0 \,\,\, \forall j \in \Upsilon\, ; \; \gamma = \gamma_0 \text{ on } \partial [1/R^2,1]^l \Big\}
\end{gather*}
and
$$
	 b_{ \lambda, \Upsilon } = \inf_{ \gamma \in \Gamma_\ast } \max_{ \textbf{t}\in [1/R^2,1]^l } \phi_\lambda \big( \gamma ( \textbf{t} ) \big).
$$

\vspace{.3cm}
Next, our intention is proving an important relation among $ b_{ \lambda, \Upsilon } $, $c_\Upsilon$ and $c_{\lambda, \Upsilon}$. However, to do this, we need to some technical lemmas. The arguments used are the same found in \cite{Alves}, however for reader's convenience we will repeat their proofs

\begin{lemma} \label{solution's existence}
   For all $ \gamma \in \Gamma_\ast $, there exists $ (s_1, \ldots, s_l ) \in [1/R^2,1]^l $ such that
$$
   \phi'_{ \lambda,\Upsilon } \big( \gamma ( s_1, \ldots, s_l ) \big) \big( \tilde{\gamma}_j ( s_1, \ldots, s_l ) \big) = 0, \, \forall j \in \Upsilon
$$
where
$$
\tilde{\gamma}_j( t_1, \ldots, t_l )(x)=
\left\{
\begin{array}{l}
\gamma ( t_1, \ldots, t_l )(x), \quad x \in \Omega'_j\\
0, \quad x \in \Omega'_\Upsilon \setminus \Omega'_j
\end{array}
\right.
$$
\end{lemma}

\noindent {\bf Proof}
   Given $ \gamma \in \Gamma_\ast $, consider $ \widehat{\gamma} \colon [1/R^2,1]^l \to \mathbb R^l $ such that
$$
   \widehat{\gamma} ( \textbf{t} ) = \Big( \phi'_{ \lambda,\Upsilon } \big( {\gamma}( \textbf{t} ) \big) \tilde
	{\gamma}_1 ( \textbf{t} ), \ldots, \phi'_{ \lambda,\Upsilon } \big( \gamma ( \textbf{t} ) \big) \tilde{\gamma}_l ( \textbf{t} ) \Big), \text{ where } \textbf{t} = ( t_1, \ldots, t_l ).
$$
For $ \textbf{t} \in \partial [1/R^2,1]^l $, it holds 
\begin{equation} \label{R0} 
\widehat{\gamma} ( \textbf{t} ) = \widehat{\gamma_0} ( \textbf{t} ), 
\end{equation}
where
$$
 \widehat{\gamma_0} ( \textbf{t} ) = \Big( I'_{ \Upsilon } \big( {\gamma_0}( \textbf{t} ) \big) t_1Rw_1, \ldots, I'_{\Upsilon } \big( \gamma_0 ( \textbf{t} ) \big) t_lRw_l \Big). 
$$
Now, lemma follows using (\ref{R}), (\ref{R0}) and Miranda's Theorem \cite{Miranda}.
\qed

\begin{proposition} \label{blambdagamma}
  \mbox{}
\begin{enumerate}
   \item[(i)] $ c_{ \lambda,\Upsilon } \le b_{ \lambda,\Upsilon } \le c_\Upsilon, \, \forall \lambda \ge 1 $;
	 \item[(ii)] $ b_{ \lambda,\Upsilon } \to c_\Upsilon, \text{ as } \lambda \to \infty $;
	 \item[(iii)] $ \phi_\lambda \big( \gamma({\bf{t}}) \big) < c_\Upsilon, \, \forall \lambda \ge 1, \gamma \in \Gamma_\ast \text{ and } {\bf{t}} = (t_1, \ldots, t_l ) \in \partial [1/R^2,1]^l $.
\end{enumerate}	
\end{proposition}

\noindent {\bf Proof}
\begin{enumerate}
   \item[(i)] Since $ \gamma_0 \in \Gamma_\ast $, by (\ref{C01}), 
	       $$
				    b_{ \lambda,\Upsilon } \le \max_{ ( t_1, \ldots, t_l ) \in [1/R^2,1]^l } \phi_\lambda ( \gamma_0 ( t_1, \ldots, t_l ) ) \leq  \max_{ ( t_1, \ldots, t_l )\in [1/R^2,1]^l} I_\Upsilon ( \sum_{j=1}^{l}t_j R w_j ) = c_\Upsilon.
				 $$
				 Now, fixing $ {\bf s} = (s_1, \ldots, s_l) \in [1/R^2,1]^l $ given in Lemma \ref{solution's existence} and recalling that
				 $$
				    c_{ \lambda,\Upsilon } =\inf_{u \in \mathcal M'_{\Upsilon}}\phi_{\lambda,\Upsilon}(u)
$$
where
$$
\mathcal M'_{\Upsilon}=\{u\in \mathcal N'_{\Upsilon}: \phi_{\lambda,\Upsilon}'(u)u_j=0 \mbox { and }
u_{j}\neq 0 \,\,\, \forall j \in \Upsilon \},
$$
				 it follows that
				 $$
				    \phi_{ \lambda,\Upsilon } ( \gamma( {\bf s } )) \geq c_{ \lambda,\Upsilon }.
				 $$
			From \eqref{til F estimate},
				 $$
				    \phi_{ \lambda, \mathbb R^3 \setminus \Omega'_\Upsilon } (u) \geq 0, \, \forall u \in H^{1} \big( \mathbb R^3 \setminus \Omega'_\Upsilon \big),
				 $$
which leads to
				 $$
				    \phi_\lambda \big( \gamma( {\textbf t} ) \big) \geq \phi_{ \lambda,\Upsilon } ( \gamma( {\textbf t} )) , \, \forall \textbf{t} = (t_1, \ldots, t_l) \in [1/R^2,1]^l.
				 $$
				 Thus
				 $$
				    \max_{ ( t_1, \ldots, t_l )\in [1/R^2,1]^l } \phi_\lambda \big( \gamma ( t_1, \ldots, t_l ) \big) \geq \phi_{\lambda,\Upsilon} \big( \gamma( \textbf{s} ) \big) \ge c_{ \lambda,\Upsilon },
				 $$
			showing that
				 $$
				    b_{ \lambda,\Upsilon } \ge c_{ \lambda,\Upsilon }.
				 $$
	 \item[(ii)] This limit is clear by the previous items, since we already know $ c_{ \lambda,\Upsilon } \to c_\Upsilon $, as $ \lambda \to \infty $;
	 \item[(iii)] For $ \textbf{t} = ( t_1, \ldots, t_l ) \in \partial [1/R^2,1]^l $, it holds $ \gamma ( \textbf{t} ) = \gamma_0 ( \textbf{t} ) $. From this,
	       $$
				    \phi_\lambda \big( \gamma ( \textbf{t} ) \big) = I_\Upsilon ( \gamma_0 ( \textbf{t} )).
				 $$
				 From (\ref{C02}) and (\ref{R}), 
				 $$
				   \phi_\lambda \big( \gamma ( \textbf{t} ) \big) \le c_\Upsilon - \epsilon,
				 $$
				 for some $ \epsilon > 0 $, so (iii) holds.
\end{enumerate}
\qed


\section{Proof of the main theorem}


To prove Theorem \ref{main}, we need to find nonnegative solutions $ u_\lambda $ for large values of $ \lambda $, which converges to a least energy solution of $(P)_{\infty,\Upsilon}$ as $ \lambda \to \infty $. To this end, we will show two propositions which together with the Propositions \ref{(PS) infty condition} and  \ref{P:boundedness of the solutions} will imply that Theorem \ref{main} holds.

Henceforth, we denote by
$$
\Theta=\left\{u \in E_\lambda\,:\, \|u\|_{\lambda, \Omega'_j}> \frac{\tau}{8R} \,\,\, \forall j \in \Upsilon \right\}
$$
and
$$
	\phi_\lambda^{ c_\Upsilon } = \big\{ u \in E_\lambda \, ; \, \phi_\lambda(u) \le c_{ \Upsilon } \big\}.
$$

Moreover, we fix $\delta=\frac{\tau}{48R}$, $ \mu >0 $ and 
\begin{equation} \label{A}
   {\cal A}_\mu^\lambda = \left\{ u \in \Theta_{2\delta} \,:\, | \phi_{ \lambda}(u)-c_\Upsilon | \leq \mu\right\}.
\end{equation}
We observe that
$$
   w_\Upsilon \in {\cal A}_\mu^\lambda \cap \phi_\lambda^{ c_\Upsilon },
$$
showing that $ {\cal A}_\mu^\lambda \cap \phi_\lambda^{ c_\Upsilon } \ne \emptyset $.

Our next result shows an important uniform estimate of $ \big\| \phi'_{ \lambda }(u) \big\|_{E^{*}_{\lambda}} $ on the region $ \left( {\cal A}_{ 2 \mu }^\lambda \setminus {\cal A}_\mu^\lambda \right) \cap \phi_\lambda^{ c_\Upsilon } $.

\begin{proposition} \label{derivative estimate}
For each  $ \mu > 0 $, there exist $ \Lambda_\ast \ge 1 $ and $ \sigma_0 >0   $ independent of $ \lambda $ such that
\begin{equation}
   \big\| \phi'_{ \lambda }(u) \big\|_{E^{*}_{\lambda}} \ge \sigma_0, \text{ for } \lambda \ge \Lambda_\ast \text{ and all } u \in \left( {\cal A}_{ 2 \mu }^\lambda \setminus {\cal A}_\mu^\lambda \right) \cap \phi_\lambda^{ c_\Upsilon }.
\end{equation}
\end{proposition}

\noindent {\bf Proof}
   We assume that there exist $ \lambda_n \to \infty $ and $ u_n \in \left( {\cal A}_{ 2 \mu }^{\lambda_n} \setminus {\cal A}_\mu^{\lambda_n} \right) \cap \phi_{\lambda_n}^{ c_\Upsilon } $ such that
$$
   \big\| \phi'_{ \lambda_n }(u_n) \big\|_{E^{*}_{\lambda_n}} \to 0.
$$
Since $ u_n \in {\cal A}_{ 2 \mu }^{ \lambda_n } $, this implies $( \|u_n\|_{ \lambda_n }) $ is a bounded sequence and, consequently, it follows that $ \big( \phi_{ \lambda_n }(u_n) \big) $ is also bounded. Thus, passing a subsequence if necessary, we can assume that $(\phi_{ \lambda_n }(u_n)) $ converges. Thus, from Proposition \ref{(PS) infty condition}, there exists $ 0 \le u \in H^{1}_0( \Omega_\Upsilon ) $ such that $u$ is a solution for $ (P)_\Upsilon $,
$$
u_n \to u \,\, \text{in} \,\, H^{1}(\mathbb{R}^{3}), \,\,	 \|u_n\|_{ \lambda_n, \mathbb R^3 \setminus \Omega_\Upsilon } \to 0 \, \text{ and } \,  \phi_{ \lambda_n} (u_n) \to I_\Upsilon(u).
$$
Recalling that $(u_n) \subset \Theta_{2 \delta}$, we derive that
$$
\|u_n\|_{\lambda_n,\Omega'_j} > \frac{\tau}{12R} \,\,\, \forall j \in \Upsilon.
$$
Then, taking the limit of $n \to +\infty$, we find
$$
\|u\|_j \geq \frac{\tau}{12R} \,\,\, \forall j \in \Upsilon,
$$
yields  $ u_{ |_{ \Omega_j } } \ne 0 $ for all $j \in \Upsilon$ and $I'_\Upsilon(u)=0 $. Consequently, by (\ref{tau0}),
$$
\|u\|_{j} > \frac{\tau}{8R} \,\,\, \forall j \in \Upsilon.
$$
This way, $I_\Upsilon (u) \geq c_{\Upsilon}$. But since $ \phi_{ \lambda_n } (u_n) \leq c_\Upsilon $ and $\phi_{ \lambda_n } (u_n) \to I_\Upsilon(u)$,  for $ n $ large, it holds
$$
\|u_n\|_j > \frac{\tau}{8R} \,\,\,  \text{ and } \, \left| \phi_{ \lambda_n}(u_n)-c_\Upsilon \right| \leq \mu, \, \forall j \in \Upsilon.
$$
So $ u_n \in {\cal A}_\mu^{\lambda_n} $, obtaining a contradiction. Thus, we have completed the proof. \qed

In the sequel, $\mu_1,\mu^*$ denote the following numbers
$$
\min_{{\bf t } \in \partial [1/R^2, 1]^l}|I_{\Upsilon}(\gamma_0 ({\bf t }))-c_\Upsilon|=\mu_{1}>0
$$
and
$$
\mu^*=\min\{\mu_1, \delta, {r}/{2}\},
$$
where $\delta$ were given (\ref{A}) and
$$
  r = R^{ 2} \left( \frac{1}{2}-\frac{1}{\theta} \right)^{-1} c_\Upsilon.
$$
Moreover, for each $s>0$, ${ B}_{s}^\lambda$ denotes the set
$$
{ B}_{s}^\lambda = \big\{ u \in E_\lambda \, ; \, \|u\|_{\lambda} \leq s \big\} \,\,\, \text{for} \,\,\, s>0.
$$

\begin{proposition} \label{P}
Let $\mu > 0 $ small enough and $ \Lambda_\ast \ge 1 $ given in the previous proposition. Then, for $ \lambda \ge \Lambda_\ast $, there exists a solution $ u_\lambda $ of $ (A_\lambda) $ such that $ u_\lambda \in {\cal A}_\mu^\lambda \cap \phi_\lambda^{ c_\Upsilon } \cap { B}_{r+1}^\lambda$.

\end{proposition}

\noindent {\bf Proof}
  Let $ \lambda \ge \Lambda_\ast $. Assume that there are no critical points of $ \phi_\lambda $ in $ {\cal A}_\mu^\lambda \cap \phi_\lambda^{ c_\Upsilon } \cap { B}_{r+1}^\lambda $. Since $ \phi_\lambda $ verifies the $ (PS) $ condition, there exists a constant $ d_\lambda > 0 $ such that
$$
   \big\| \phi'_\lambda(u) \big\|_{E^{*}_{\lambda}} \ge d_\lambda, \text{ for all } u \in {\cal A}_\mu^\lambda \cap \phi_\lambda^{ c_\Upsilon } \cap { B}_{r+1}^\lambda.
$$
From Proposition \ref{derivative estimate}, 
$$
   \big\| \phi'_\lambda(u) \big\|_{E^{*}_{\lambda}} \ge \sigma_0, \text{ for all } u \in \left( {\cal A}_{ 2 \mu }^{\lambda} \setminus {\cal A}_\mu^{\lambda} \right) \cap \phi_{\lambda}^{ c_\Upsilon },
$$
where $ \sigma_0 > 0 $ does not depend on $ \lambda $. In what follows,  $ \Psi \colon E_\lambda \to \mathbb R $ is a continuous functional verifying
$$
   \Psi(u) =  1, \text{ for } u \in {\cal A}_{\frac{3}{2} \mu}^\lambda \cap \Theta_\delta \cap B^{\lambda}_{r},
$$
$$
\ \Psi(u) = 0, \text{ for } u \notin {\cal A}_{2 \mu}^\lambda \cap \Theta_{2\delta} \cap B^{\lambda}_{r+1}
$$
and
$$
0 \le \Psi(u) \le 1, \, \forall u \in E_\lambda.
$$
We also consider $ H \colon \phi_\lambda^{ c_\Upsilon } \to E_\lambda $ given by
$$
   H(u) =
\begin{cases}
   - \Psi(u) \big\| Y(u) \big\|^{ -1 } Y(u), \text{ for } u \in {\cal A}_{2 \mu}^\lambda \cap B^{\lambda}_{r+1}, \\
   \phantom{- \Psi(u) \big\| Y(u) \big\|^{ -1 } Y()} 0, \text{ for } u \notin {\cal A}_{2 \mu}^\lambda \cap B^{\lambda}_{r+1}, \\
\end{cases}
$$
where $ Y $ is a pseudo-gradient vector field for $ \Phi_\lambda $ on $ {\cal K} = \left\{ u \in E_\lambda \, ; \, \phi'_\lambda(u) \ne 0 \right\} $. Observe that  $ H $ is well defined, once $ \phi'_\lambda(u) \ne 0 $, for $ u \in {\cal A}_{2 \mu}^\lambda \cap \phi_\lambda^{ c_\Upsilon } $. The inequality
$$
   \big\| H(u) \big\| \le 1, \, \forall \lambda \ge \Lambda_* \text{ and } u \in \phi_\lambda^{ c_\Upsilon },
$$
guarantees that the deformation flow $ \eta \colon [0, \infty) \times \phi_\lambda^{ c_\Upsilon } \to \phi_\lambda^{ c_\Upsilon } $ defined by
$$
   \frac{d \eta}{dt} = H(\eta), \ \eta(0,u) = u \in \phi_\lambda^{ c_\Upsilon }
$$
verifies
\begin{gather}
   \frac{d}{dt} \phi_\lambda \big( \eta(t,u) \big) \le - \frac{1}{2} \Psi \big( \eta(t,u) \big) \big\| \phi'_\lambda \big( \eta(t,u) \big) \big\| \le 0, \label{eta derivative}\\
	 \left\| \frac{d \eta}{dt} \right\|_\lambda = \big\| H(\eta) \big\|_\lambda \le 1
\end{gather}
and
\begin{equation} \label{eta}
   \eta(t,u) = u \text{ for all } t \ge 0 \text{ and } u \in \phi_\lambda^{ c_\Upsilon } \setminus {\cal A}_{2 \mu}^\lambda \cap B^{\lambda}_{r+1}.
\end{equation}

Next, we study two paths, which are relevant for what follows: \\

$ \noindent \bullet $ The path $ {\bf t} \mapsto \eta \big( t, \gamma_0( {\bf t} ) \big), \text{ where } \textbf{t} = (t_1,\ldots,t_l) \in [1/R^2, 1]^l $.

\vspace{0.5 cm}

Thereby, if $\mu \in (0,\mu^*)$, we have that
$$
   \gamma_0( {\bf t } ) \notin {\cal A}_{2 \mu}^\lambda, \, \forall {\bf t } \in \partial [1/R^2, 1]^l.
$$
Since
$$
   \phi_\lambda \big( \gamma_0( {\bf t } ) \big) < c_\Upsilon, \, \forall {\bf t } \in \partial [1/R^2, 1]^l,
$$
from (\ref{eta}), it follows that
$$
   \eta \big( t, \gamma_0( {\bf t} ) \big) = \gamma_0( {\bf t} ), \, \forall {\bf t} \in \partial [1/R^2, 1]^l.
$$
So, $ \eta \big( t, \gamma_0( {\bf t} ) \big) \in \Gamma_\ast $, for each $ t \ge 0 $.

\vspace{0.5 cm}

$ \noindent \bullet $ The path $ {\bf t} \mapsto \gamma_0( {\bf t} ), \text{ where } \textbf{t} = (t_1,\ldots,t_l) \in [1/R^2, 1]^l $.

\vspace{0.5 cm}

We observe that
$$
   \text{supp} \big( \gamma_0 ( {\bf t} ) \big)\subset \overline{\Omega_\Upsilon}
$$
and
$$
   \phi_\lambda \big( \gamma_0 ( {\bf t} ) \big) \text{ does not depend on }  \lambda \ge 1,
$$
for all  $ {\bf t} \in [1/R^2, 1]^l $. Moreover,
$$
   \phi_\lambda \big( \gamma_0 ( {\bf t} ) \big) \le c_\Upsilon, \, \forall {\bf t} \in [1/R^2, 1]^l
$$
and
$$
   \phi_\lambda \big( \gamma_0 ( {\bf t} ) \big) = c_\Upsilon \text{ if, and only if, } t_j = \frac{1}{R}, \, \forall j \in \Upsilon.
$$
Therefore
$$
   m_0 = \sup \left\{ \phi_\lambda(u) \, ; \, u \in \gamma_0 \big( [1/R^2,1]^l \big) \setminus A_\mu^\lambda \right\}
$$
is independent of $ \lambda $ and $ m_0 < c_\Upsilon $. Now, observing that there exists $ K_\ast > 0 $ such that
$$
   \big| \phi_{ \lambda }(u) - \phi_{ \lambda}(v) \big| \le K_* \| u-v \|_{ \lambda }, \, \forall u,v \in {\cal B}_r^\lambda ,
$$
we derive
\begin{equation} \label{max estimate}
    \max_{ {\bf t } \in [1/R^2,1]^l } \phi_\lambda \Big( \eta \big( T, \gamma_0 ( {\bf t} ) \big) \Big) \le \max \left\{ m_0, c_\Upsilon-\frac{1}{2 K_\ast} \sigma_0 \mu \right\},
\end{equation}
for $ T > 0 $ large.

In fact, writing $ u = \gamma_0( {\bf t} ) $, $ {\bf t } \in [1/R^2,1]^l $, if $ u \notin A_\mu^\lambda $, from (\ref{eta derivative}),
$$
   \phi_\lambda \big( \eta( t, u ) \big) \le \phi_\lambda (u) \le m_0, \, \forall t \ge 0,
$$
and we have nothing more to do. We assume then $ u \in A_\mu^\lambda $ and set
$$
   \widetilde{\eta}(t) = \eta (t,u), \ \widetilde{d_\lambda} = \min \left\{ d_\lambda, \sigma_0 \right\} \text{ and } T = \frac{\sigma_0 \mu}{K_\ast \widetilde{d_\lambda}}.
$$
Now, we will analyze the ensuing cases: \\

\noindent {\bf Case 1:} $ \widetilde{\eta}(t) \in {\cal A}_{\frac{3}{2} \mu}^\lambda \cap \Theta_\delta \cap B^{\lambda}_{r}, \, \forall t \in [0,T] $.

\noindent {\bf Case 2:} $ \widetilde{\eta}(t_0) \notin {\cal A}_{\frac{3}{2} \mu}^\lambda \cap \Theta_\delta \cap B^{\lambda}_{r}, \text{ for some } t_0 \in [0,T] $. \\

\noindent {\bf Analysis of  Case 1}

In this case, we have $ \Psi \big( \widetilde{\eta}(t) \big) = 1 $ and $ \big\| \phi'_\lambda \big( \widetilde{\eta}(t) \big) \big\| \ge \widetilde{d_\lambda} $ for all $ t \in [0,T] $. Hence, from (\ref{eta derivative}),
$$
   \phi_\lambda \big( \widetilde{\eta}(T) \big) = \phi_\lambda (u) + \int_0^T \frac{d}{ds} \phi_\lambda \big( \widetilde{\eta}(s) \big) \, ds \le c_\Upsilon - \frac{1}{2} \int_0^T \widetilde{d_\lambda} \, ds,
$$
that is,
$$
   \phi_\lambda \big( \widetilde{\eta}(T) \big) \le c_\Upsilon - \frac{1}{2} \widetilde{d_\lambda} T = c_\Upsilon - \frac{1}{2 K_\ast} \sigma_0 \mu,
$$
showing (\ref{max estimate}). \\

\noindent {\bf Analysis of Case 2}:  In this case we have the following situations:
\\

\noindent {\bf (a)}: There exists $t_2 \in [0,T]$ such that $\tilde{\eta}(t_2) \notin \Theta_\delta$, and thus, for $t_1=0$ it follows that
$$
\|\tilde{\eta}(t_2)-\tilde{\eta}(t_1)\| \geq \delta > \mu,
$$
because $\tilde{\eta}(t_1)=u \in \Theta$. \\

\vspace{0.5 cm}

\noindent {\bf (b)}: There exists $t_2 \in [0,T]$ such that $\tilde{\eta}(t_2) \notin B^{\lambda}_r$, so that for $t_1=0$, we get
$$
\|\tilde{\eta}(t_2)-\tilde{\eta}(t_1)\| \geq r > \mu,
$$
because $\tilde{\eta}(t_1)=u \in B^{\lambda}_r$. \\

\noindent {\bf (c)}: \, $\tilde{\eta}(t) \in \Theta_\delta \cap B^{\lambda}_r$ for all $t \in [0,T]$, and there are $0 \leq t_1 \leq t_2 \leq T$ such that $\tilde{\eta}(t) \in  {\cal A}_{\frac{3}{2} \mu}^\lambda \setminus {\cal A}_\mu^\lambda$ for all $t \in [t_1,t_2]$ with
$$
|\phi_\lambda(\tilde{\eta}(t_1))-c_\Upsilon|=\mu \,\,\, \mbox{and} \,\,\, |\phi_\lambda(\tilde{\eta}(t_2))-c_\Upsilon|=\frac{3\mu}{2}
$$
From definition of $K_\ast$, we have
$$
   \| w_2-w_1 \| \ge \frac{1}{K_\ast} \big| \phi_{ \lambda} (w_2) - \phi_{ \lambda} (w_1) \big| \ge \frac{1}{2 K_\ast} \mu.
$$
Then, by  mean value theorem, $ t_2-t_1 \ge \frac{1}{2 K_\ast} \mu $ and, this form,
$$
   \phi_\lambda \big( \widetilde{\eta}(T) \big) \le \phi_\lambda(u) - \int_0^T \Psi \big( \widetilde{\eta}(s) \big) \big\| \phi'_\lambda \big( \widetilde{\eta}(s) \big) \big\| \, ds
$$
implying
$$
   \phi_\lambda \big( \widetilde{\eta}(T) \big) \le c_\Upsilon - \int_{t_1}^{t_2} \sigma_0 \, ds = c_\Upsilon - \sigma_0 (t_2-t_1) \le c_\Upsilon - \frac{1}{2 K_\ast} \sigma_0 \mu,
$$
which proves (\ref{max estimate}). Fixing $ \widehat{\eta} (t_1, \ldots, t_l) = \eta \big( T, \gamma_0 (t_1,\ldots,t_l) \big) $, we have that
$\widehat{\eta}(t_1, \ldots, t_l) \in \Theta_{2\delta}$, and so, $\widehat{\eta}(t_1, \ldots, t_l)|_{\Omega'_j} \not= 0$ for all $j \in \Upsilon$. Thus, $ \widehat{\eta} \in \Gamma_\ast $, leading to
$$
   b_{ \lambda, \Gamma } \le \max_{ (t_1,\ldots,t_l) \in [1/R^2, 1] } \phi_\lambda \big( \widehat{\eta} (t_1,\ldots,t_l) \big) \le \max \left\{ m_0, c_\Upsilon - \frac{1}{2 K_\ast} \sigma_0 \mu \right\} < c_\Upsilon,
$$
which contradicts the fact that $ b_{ \lambda, \Upsilon } \to c_\Upsilon $.
\qed

\vspace{.5cm}
\noindent {\bf [Proof of Theorem \ref{main}]}
According Proposition \ref{P}, for $\mu \in (0, \mu^*)$ and $ \Lambda_\ast \ge 1 $, there exists a solution $ u_\lambda $ for $ (A_\lambda) $ such that $ u_\lambda \in {\cal A}_\mu^\lambda \cap \phi_\lambda^{ c_\Upsilon } $, for all $\lambda \geq \Lambda_*$. \\

\noindent {\bf Claim:}
There are $\lambda_0 \geq \Lambda_*$ and $\mu_0>0$ small enough, such that $u_\lambda$ is a solution for $ (P)_\lambda $ for $\lambda \geq \Lambda_0$ and $\mu \in (0, \mu_0)$.

Indeed, fixed $\mu \in (0, \mu_0)$, assume by contradiction that there are  $ \lambda_n \to \infty $, such that $(u_{\lambda_n})$ is not a solution for $(P)_{\lambda_n}$. From Proposition \ref{P}, the sequence $ (u_{\lambda_n}) $ verifies:
\begin{enumerate}
   \item[(a)] $ \phi'_{ \lambda_n }(u_{\lambda_n}) = 0, \, \forall n \in \Bbb N $;
	 \item[(b)] $ \|u_n\|^{2}_{\lambda_n,  \Bbb R^3 \setminus \Omega_\Upsilon }(u_{\lambda_n}) \to 0$;
	 \item[(c)] $ \phi_{ \lambda_n } (u_{\lambda_n}) \to d \leq c_\Upsilon. $
\end{enumerate}	
The item (b) ensures we can use Proposition \ref{P:boundedness of the solutions} to deduce $ u_{\lambda_n} $ is a solution for
$(P)_{\lambda_n} $, for large values of $ n $, which is a contradiction, showing this way the claim. \\

Now, our goal is to prove the second part of the theorem. To this end, let  $(u_{\lambda_n})$ be a sequence verifying the above limits. A direct computation gives $ \phi_{ \lambda_n }(u_{ \lambda_n } ) \to d$ with $d \leq c_\Upsilon $. This way, using Proposition \ref{(PS) infty condition} combined with item (c), we derive $( u_{ \lambda_n } )$ converges in $ H^{1}(\Bbb R^3) $ to a function $ u \in H^{1}(\mathbb{R}^3) $, which satisfies $ u = 0 $ outside $ \Omega_\Upsilon $ and $ u_{|_{\Omega_j}} \not= 0, \, j \in \Upsilon $, and $u$ is a positive solution for
         $$
				 \begin{cases}
		        M\biggl(\displaystyle\int_{\Omega_\Upsilon}(|\nabla u|^{2}+u^{2}) dx\biggl)(- \Delta u + u )= f(u), \text{ in } \Omega_\Upsilon, \\
		        u \in H^{1}_0 ( \Omega_\Upsilon ),
		       \end{cases}
\eqno{(P)_{\infty, \Upsilon}}
$$
and so,
$$
I_{\Upsilon}(u) \geq c_{\Upsilon}.
$$
On the other hand, we also know that
$$
\phi_{ \lambda_n }(u_{ \lambda_n } ) \to I_{\Upsilon}(u),
$$
implying that
$$
I_{\Upsilon}(u)=d \,\,\, \mbox{and} \,\,\, d \geq c_{\Upsilon}.
$$
Since $d \leq c_{\Upsilon}$, we deduce that
$$
I_{\Upsilon}(u)=c_{\Upsilon},
$$
showing that $u$ is a least energy solution for $(P)_{\infty, \Upsilon}$. Consequently, $u$ is a least energy solution for the problem	
	$$
      	    \begin{cases}
		     M\biggl(\displaystyle\int_{\Omega_\Upsilon}(|\nabla u|^{2}+u^{2}) dx\biggl)(- \Delta u + u)  = f(u), \text{ in } \Omega_\Upsilon, \\
						     u \in H^{1}_0 ( \Omega_\Upsilon ).
			\end{cases}
  	$$
\qed


\vspace{0.2cm}

\begin{thebibliography}{99}

































 \bibitem{Alves} C.O. Alves,{\it  \, Existence of multi-bump solutions for a class of quasilinear problems,} { Adv. Nonlinear Stud.} {6}(2006), 491-509 .

\bibitem{alvescorrea} C.O. Alves and  F.J.S.A. Corr\^{e}a , {\it On
existence of solutions for a class of problem involving a nonlinear
operator}, Comm. Appl. Nonlinear Anal. 8(2001), 43-56.

\bibitem{alvescorreama}
C.O. Alves, F.J.S.A. Corr\^{e}a  and T.F. Ma,  {\it Positive
solutions for a quasilinear elliptic equation of Kirchhoff type,}
{ Comput. Math. Appl. 49(2005)85-93}.

\bibitem{AlvesFigueiredo}
C. O. Alves  and G. M. Figueiredo,  {\it Nonlinear perturbations of a
periodic Kirchhoff equation in $\mathbb{R}^{N}$.} { Nonlinear
Anal. 75(2012)2750-2759}.

\bibitem{AY} C.O. Alves and M. Yang, {\it Existence of positive multi-bump solutions for a  Schr\"odinger-Poisson system in 
$\mathbb{R}^{3}$}, arXiv:1501.02930v1


\bibitem{Azzollini} A. Azzollini, {\it The elliptic Kirchhoff equation in $\mathbb{R}^{N}$ perturbed by a local nonlinearity.}
{ Differential Integral Equations 25 (2012), 543-554.}

\bibitem{bl} {H. Berestycki} and {P.L. Lions, } {\it Nonlinear scalar field equations, I - existence of a ground state },
 {Arch. Rat. Mech. Analysis} {82} {(1983)}, {313--346}.

\bibitem{Cammaroto} F. Cammaroto  and L. Vilasi, {\it On a Schr\"{o}dinger-Kirchhoff-type equation involving the
$p(x)$-Laplacian}, Nonlinear Anal. 81 (2013), 42-53.

\bibitem{Cheng} B. Cheng, X. Wu and J. Liu, {\it Multiple solutions for a class of Kirchhoff
type problems with concave nonlinearity}, NoDEA Nonlinear
Differential Equations Appl. 19 (2012), 521-537.

\bibitem{Chen} C. Chen, H. Song and  Z. Xiu, {\it Multiple solutions for p-Kirchhoff equations in
$\mathbb{R}^{N}$}, Nonlinear Anal. 86 (2013), 146-156.


\bibitem{SChen} S. Chen and L. Li, {\it Multiple solutions for the nonhomogeneous Kirchhoff equation
on $\mathbb{R}^{N}$}, Nonlinear Anal. Real World Appl. 14 (2013),
1477-1486.

\bibitem{DelPinoFelmer} M. del Pino and P.L. Felmer, Local mountain passes for semilinear elliptic problems in unbounded domains, { \ Calc. Var. PDE} {4} (1996), 121-137.


\bibitem{DingTanaka} Y.H. Ding and K. Tanaka, {\it Multiplicity of positive solutions of a nonlinear Schr\"odinger equation}, { Manuscripta Math.} {112} (2003), 109-135.


\bibitem{jmaa} G.M. Figueiredo, {\it Existence of positive solution for a Kirchhoff  problem type with
critical growth via truncation argument}, { J. Math. Anal. Appl.
401 (2013), 706-713.}

\bibitem{Cristian} G.M. Figueiredo, C. Morales, J. R. Santos Junior and A. Suarez, {\it 
Study of a nonlinear Kirchhoff equation with non-homogeneous material }, { J. Math. Anal. Appl.
416 (2014), 597-608.}

\bibitem{Giovany1} G.M. Figueiredo and J. R. Santos Junior, {\it 
Multiplicity and concentration behavior of positive solutions for a Schrodinger-Kirchhoff type problem 
via penalization method }, {ESAIM: Control, Otimiz. Calc. Variat. 20 (2014), 389-415.}

\bibitem{Giovany2} G.M. Figueiredo, N. Ikoma and J. R. Santos Junior, {\it 
Existence and concentration result for the Kirchhoff type equations with general nonlinearities }, 
{ARMA 213 (2014), 931-979}.

\bibitem{Giovany3} G.M. Figueiredo and R. G. Nascimento, {\it 
Existence of a nodal solution with minimal energy for a Kirchhoff equation }, 
{ Math. Nachrichten 288 (2015), 48-60.}

\bibitem{Gon} L. Gongbao, {\it Some properties of weak solutions of nonlinear scalar field equations}, Ann. Acad. Sci. Fenn. Math. 14 (1989), 27-36.

\bibitem{He}
X. He  and W. Zou,  {\it Existence and concentration of positive
solutions for a Kirchhoff equation in $\R^{3}$}. { J.
Differential Equations, 252(2012)1813-1834}.


\bibitem{kirchhoff}
G. Kirchhoff, {\it Mechanik, Teubner,Leipzig, 1883}.

\bibitem{SS} S. Liang and S. Shi, {\it Existence of multi-bump solutions for a class of Kirchhoff type problems in $\mathbb{R}^{3}$}, J. Math. Phys. 54, 121510 (2013); doi 10.1063/1.4850835


\bibitem{Li}
Y. Li, F. Li  and J. Shi,  {\it Existence of a positive solution to
Kirchhoff type problems without compactness conditions}.  {J.
Differential Equations 253 (2012), 2285-2294}.

\bibitem{Liao}
YJia-Feng Liao, Peng Zhang, Jiu Liu and Chun-Lei Tang {\it Existence and multiplicity of positive solutions for a class of Kirchhoff 
type problems with singularity}.  { J. Math. Anal. Appl. 430 (2015), 1124-1148}
  

\bibitem{Liu} X. Liu and Y. Sun, {\it Multiple positive solutions for Kirchhoff type problems with
singularity}, Commun. Pure Appl. Anal. 12 (2013), 721-733.

\bibitem{Liu1} Z. Liu and S. Guo, {\it Existence of positive ground state solutions for Kirchhoff type problems}, 
Nonlinear Anal. 120 (2015), 1-13.

\bibitem{ma}
T.F. Ma, {\it Remarks on an elliptic equation of Kirchhoff type}.
{ Nonlinear Anal.,  63 (2005)1967-1977}.

\bibitem{Mao} A. Mao and Z. Zhang, {\it Sign-changing and multiple solutions of Kirchhoff type problems
without the P.S. condition}, Nonlinear Anal. 70 (2009), 1275-1287.

\bibitem{Mao1} A. Mao and Shixia Luan,  {\it Sign-changing solutions of a class of nonlocal quasilinear elliptic
boundary value problems}, J. Math. Anal. Appl., 383 (2011) 239-243.


\bibitem{Miranda} {C. Miranda}, {Un' osservazione su un teorema di Brouwer,}{ Bol. Un. Mat. Ital. 3 (1940)}, 5-7.

\bibitem{Naimen}
D. Naimen, {\it On the Brezis-Nirenberg problem with a Kirchhoff type perturbation.}, Adv. Nonlinear Stud. 15 (2015), 135-156.


\bibitem{Naimen1}
D. Naimen, {\it The critical problem of Kirchhoff type elliptic equations in dimension four.}, 
J. Differential Equations, 257 (2014), 1168-1193. 

\bibitem{Shuai}
W. Shuai, {\it Sign-changing solutions for a class of Kirchhoff-type problem in bounded domains}, J. Differential Equations 254 (2015), 1256-1274. 





\bibitem{Xiang} M. Xiang, B. Zhang and X. Guo, {\it Infinitely many solutions for a fractional Kirchhoff type 
problem via Fountain Theorem}, Nonlinear Anal., 120 (2015),  299-313.


\bibitem{Wang1}
J. Wang , L. Tian  , J. Xu  and F. Zhang,  {\it Multiplicity and
concentration of positive solutions for a Kirchhoff type problem
with critical growth}. {J. Differential Equations 253 (2012), 2314-2351}.

\bibitem{Wang2} L. Wang,   {\it On a quasilinear Schr\"{o}dinger-Kirchhoff-type equation with radial potentials,} 
{Nonlinear Anal. 83 (2013), 58-68.}


\bibitem{W} M. Willem, Minimax Theorems, Birkh\"auser Boston, MA (1996).


\bibitem{Zhang} Z. Zhang and K. Perera, {\it Sign changing solutions of Kirchhoff type problems
via invariant sets of descent flow}, J. Math. Anal. Appl. 317 (2006),
456-463.



\end{thebibliography}
\end{document}